\documentclass[a4paper,USenglish,cleveref, autoref, thm-restate,numberwithinsect]{lipics-v2021}

\pdfoutput=1
\hideLIPIcs  

\usepackage{bm}
\usepackage{mathtools}
\usepackage{xspace}
\usepackage{xcolor}
\usepackage{csquotes}
\usepackage{cite}
\usepackage{placeins}

\definecolor{signalblue1}{RGB}{58,118,240}
\DeclareTextFontCommand{\emph}{\color{signalblue1}\em}

\DeclareMathOperator{\interior}{int}
\newcommand{\td}{\mathcal{T}}
\newtheorem{question}{Question}
\newcommand{\nice}{strong\xspace}

\newcommand{\floor}[1]{\left\lfloor #1 \right\rfloor}
\let\OldBox\Box
\renewcommand{\Box}{\mathbin{\OldBox}}
\usepackage{float}

\graphicspath{{./figures/}}

\title{Characterizing Fixed-Treewidth Planar Graphs}

\author{Laura Merker}{Karlsruhe Institute of Technology, Germany}{laura.merker2@kit.edu}{https://orcid.org/0000-0003-1961-4531}{}
\author{Samuel Schneider}{Karlsruhe Institute of Technology, Germany}{samuel.schneider@kit.edu}{https://orcid.org/0009-0002-9680-4048}{}
\author{Torsten Ueckerdt}{Karlsruhe Institute of Technology, Germany}{torsten.ueckerdt@kit.edu}{https://orcid.org/0000-0002-0645-9715}{}

\authorrunning{L. Merker, S. Schneider, T. Ueckerdt}

\Copyright{Laura Merker, Samuel Schneider, Torsten Ueckerdt}

\ccsdesc[500]{Mathematics of computing~Graphs and surfaces}

\keywords{planar graphs, treewidth, stacked triangulations, vertex arboricity}
\nolinenumbers

\begin{document}

\maketitle

\begin{abstract}
    We give the first inductive characterization of planar graphs of treewidth $k$ for each $k \geq 3$.
    Our definition generalizes the well-known inductive definition of planar $3$-trees and the recent definition of so-called \emph{planar quasi-$4$-trees} [Förster, 2024].
    As a key step, which may be of independent interest, we show that every planar triangulation $G$ of treewidth~$k$ admits a tree decomposition of width~$k$ where each adhesion set induces a cycle in $G$.
    As an application of planar quasi-$4$-trees, 
    we show that planar graphs of treewidth at most $4$ have \emph{vertex arboricity} at most $2$.
    That is, the vertex-set of every planar graph of treewidth at most $4$ can be partitioned into two sets that each induce a forest.
\end{abstract}

\section{Introduction}

Undoubtedly, treewidth is one of the most important parameters in algorithmic and structural graph theory.
Numerous investigations are concerned with graphs of small treewidth and their combinatorial properties.
The treewidth of a graph can be defined in several, seemingly different ways, e.\,g., via tree decompositions, chordal completions, clique-sums, or a variant of the cops-and-robber game (see for example this survey~\cite{treewidth_survey}).
A very simple and accessible definition involves the edge-maximal graphs of treewidth~$k$, so-called \emph{$k$-trees}, which are defined starting with $K_{k+1}$, and inductively adding a new vertex connected to exactly the vertices of some $k$-clique of the so-far constructed graph.

Planar graphs form a second central pillar of graph theory, both from a theoretical and practical point of view.
While a planar graph can have arbitrarily large treewidth, a widespread and worthwhile approach is to consider planar graphs of small treewidth.
As all graphs of treewidth at most~$2$ are planar, the combination of planarity and small treewidth becomes interesting starting with treewidth~$3$.
In fact, planar graphs of treewidth at most~$3$ are well understood.
These are exactly the subgraphs of \emph{planar $3$-trees} also known as \emph{stacked triangulations} or \emph{Apollonian networks}~\cite{planar_partial_3-trees_are_partial_planar_3-trees,a_note_on_planar_partial_3-trees}.
In particular it is well known that every planar $3$-tree can be obtained by starting with a triangle and then repeatedly adding a vertex to an inner face and connecting it to the vertices incident to the face (see for example \cite{triangulations_with_same_tw}).

The investigation of planar graphs of low treewidth has become even more vital and fundamental with the development of the \emph{product structure theory} in recent years.
Roughly speaking, product structure reduces a graph~$G$ of large treewidth to a planar graph of small treewidth by \enquote{factoring out} a path and a small clique.
Specifically, for small constants $k$ and $c$, one seeks to have $G$ as a subgraph of $H_k \boxtimes P \boxtimes K_c$, where $H_k$ is a planar graph of treewidth~$k$, and $P$ is a path.%
\footnote{The term $G\boxtimes H$ denotes the \emph{strong product} of the graphs $G$ and $H$ with $V(G\boxtimes H) = V(G) \times V(H)$ and $(u_1,v_1)(u_2,v_2) \in E(G \boxtimes H)$ if $u_1=u_2, v_1=v_2$ or $u_1u_2 \in E(G)$ and $v_1v_2\in E(H)$.}
Several well-structured graph classes admit product structure with the treewidth~$k$ of the planar factor between~$3$ and~$6$~\cite{PlanarGraphsProductStructure,ImprovedPlanarGraphProductStructureTheorem,h_framed_product_structure,ProductStructurekPlanarGraphs}.
This has been used to achieve a plethora of results for these and related graph classes~\cite{h_framed_product_structure,queue_number_42,ShorterLabelingSchemesforPlanarGraphs,ReducedBandwidth,VertexRankingPlanarGraphs,clustered_colorings_product_structure,centered_colorings_product_structure,AdjacencyLabellingforPlanarGraphsAndBeyond,Dujmovi__2020,clustered_product_structure,clustered_hadwiger_product_structure,PlanarGraphsProductStructure,OddColouringsofGraphProducts,box_dimension_product_structure,SparseUniversalGraphsForPlanarity,goetze2025strongoddcoloringminorclosed,joret2024neighborhood,twin_width_product_structure}.

\begin{figure}
    \centering
    \includegraphics{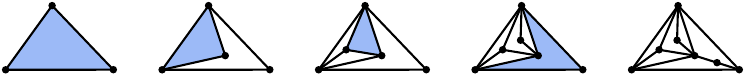}
    \caption{From left to right an example construction of a planar quasi-$4$-tree, where in each step a new vertex is inserted into the blue face of the previous step. Note that in some steps, not all possible edges are added to the new vertex, creating faces of size at most~$4$.}
    \label{fig:example_construction}
\end{figure}

\subparagraph{Our Contribution.}
While planar graphs of treewidth~$3$ are well understood, there is no clear understanding of planar graphs of treewidth~$k > 3$.
To the best of our knowledge, there is no inductive characterization of planar graphs of treewidth~$k$ for any $k > 3$.
Recently, Förster~\cite{planar_quasi_4_trees}\footnote{This paper has been withdrawn for reasons not regarding the definition of planar quasi-$4$-trees.} defined \emph{planar quasi-$4$-trees} by generalizing the inductive construction of planar $3$-trees.
Förster does this by allowing to omit some of the edges when inserting a new vertex so that faces of size at most 4 may occur, see for example \cref{fig:example_construction}.
Förster conjectured that in fact every planar graph of treewidth at most~$4$ is a subgraph of a planar quasi-$4$-tree.
As our main contribution, we generalize this construction to all $k \geq 3$ and show that a planar graph has treewidth at most~$k$ if and only if it is a subgraph of a \emph{planar quasi-$k$-tree}.
In particular, this confirms the conjecture by Förster.%
\footnote{
    Our definition of planar quasi-$k$-trees yields a slightly different construction for $k=4$ than the planar quasi-$4$-trees as defined by Förster.
    Förster's definition is equivalent to our definition of \emph{\nice} planar quasi-$4$-trees (see \cref{sec:niceness}).
    We relate these two definitions in \cref{thm:nice_4-trees}.
}
To the best of our knowledge, no such characterization of planar graphs of treewidth~$k$ was known prior to this work.
However, Seymour and Thomas~\cite{ratcatcher} and Dorn, Penninkx, Bodlaender, and Fomin~\cite{sphere_cut} give a characterization of planar graphs of \emph{branchwidth}~$b$ via \emph{sphere-cut-decompositions}.
As every graph of branchwidth~$b$ has treewidth at least~$b-1$ and at most~$\floor{\frac{3}{2}b}-1$~\cite{branchwidth}, this gives a sort of approximate characterization for planar graphs of fixed treewidth.
We discuss this in more detail in \cref{sec:comparing_sphere_cut}.

To indicate the usefulness of our construction, our second contribution is an application of the developed inductive structure of planar quasi-$4$-trees.
We look at the \emph{vertex arboricity}, also known as \emph{point arboricity}, of planar graphs of treewidth at most~$4$.
The vertex arboricity of a graph~$G$ is the smallest integer $t$ such that $V(G)$ can be partitioned into~$t$ sets such that each set induces a forest.
Chartrand, Kronk, and Wall~\cite{chartrandPointarboricityGraph1968} introduce this parameter and show that every planar graph has vertex arboricity at most~$3$.
Since then, it was shown that various subclasses of planar graphs have vertex arboricity at most~$2$, such as triangle-free planar graphs \cite{arboricity_triangle_free}, planar graphs without $4$-cycles \cite{RASPAUD20081064}, and planar graphs of degeneracy at most~$3$ \cite{arboricity_degeneracy}, which include planar 3-trees.
For general planar graphs, however, Chartrand and Kronk~\cite{point_arboricity_3} show that~$3$ is indeed tight by showing that the dual of the Tutte graph~\cite{tutte_graph} has vertex arboricity~$3$.
Since the dual of the Tutte graph has treewidth~$5$, this leaves the question open whether planar graphs of treewidth at most~$4$ have vertex arboricity at most~$2$.
We utilize planar quasi-$4$-trees to show that this is indeed the case.

\subparagraph{Organization of the Paper and Results.}\label{se:results}

We start by defining basic concepts and notation in \cref{se:preliminaries} and then provide our definition of \emph{planar quasi-$k$-trees} in \cref{sec:quasi_trees} (cf.\ \cref{def:quasi_k-trees}), which is devoted to our main theorem.

\begin{restatable}{theorem}{quasiktrees}\label{thm:quasi_k_trees}
    For every $k \geq 3$, a planar graph $G$ has treewidth at most $k$ if and only if $G$ is a subgraph of a planar quasi-$k$-tree.
\end{restatable}

We actually show the slightly stronger statement that the planar quasi-$k$-tree which is a supergraph of $G$ has the same vertex-set as $G$.
Our proof of \cref{thm:quasi_k_trees} relies on tree decompositions in which every adhesion set induces a cycle, which we call \emph{cycle-inducing} tree decompositions.
We believe that cycle-inducing tree decompositions might be of independent interest.

\begin{restatable}{theorem}{cycleinducing} \label{thm:nice_edges} \label{le:tree_decomposition_nice_edges}
    Let $G$ be a planar triangulation of treewidth $k$. 
    Then, $G$ admits a cycle-inducing tree decomposition of width $k$.
\end{restatable}

For $ k = 4 $, we strengthen \cref{thm:quasi_k_trees} in \cref{se:quasi_4-trees} by introducing \emph{\nice} planar quasi-$k$-trees (cf.\ \cref{def:strong_quasi_k-trees}) and show that for $k=4$ these yield the same maximal graphs as planar quasi-$4$-trees.

\begin{restatable}{theorem}{nicefourtrees}\label{thm:nice_4-trees}
    A planar graph has treewidth at most $4$ if and only if it is a subgraph of a \nice planar quasi-$4$-tree.
\end{restatable}

We then exploit this stronger notion to show that every planar graph of treewidth~$4$ has vertex arboricity~$2$, closing the gap between all planar graphs of treewidth 3 having vertex arboricity 2~\cite{arboricity_degeneracy} and some planar graphs of treewidth 5 having vertex arboricity 3~\cite{point_arboricity_3}.

\begin{theorem}\label{thm:arboricity}
    Let $G$ be a planar graph of treewidth at most $4$. Then, $V(G)$ can be partitioned into two sets $A$ and $B$ such that $G[A]$ and $G[B]$ are forests.
\end{theorem}

Having understood the usage of cycle-inducing tree decompositions in \cref{sec:quasi_trees,se:quasi_4-trees}, we show their existence in \cref{se:cycle_inducing_tree_decompositions}, which may be of independent interest.
Then, to complement the positive results, we show in \cref{sec:niceness} that \cref{thm:nice_4-trees} cannot be generalized to $ k > 4 $.

\begin{restatable}{theorem}{notnicetrees}\label{thm:not_nice_k_tree}
    For every $k \geq 5$, there is a planar graph of treewidth $5$ that is not a subgraph of any \nice planar quasi-$k$-tree.
\end{restatable}

\Cref{thm:nice_4-trees,thm:not_nice_k_tree} can be interpreted so that for $ k = 4 $, the degree of the vertices that are inserted in \cref{fig:example_construction} is at least 2 at the time of insertion, whereas for $ k > 4 $, degree 1 may be necessary.

We conclude this paper with some open questions and a conjecture in \cref{se:open_questions}.
Most importantly, we ask whether our inductive definition of planar quasi-$k$-trees can be used to efficiently generate planar graphs of treewidth $ k $, similarly to the inductive definition of planar 3-trees, which is used in \texttt{plantri} by Brinkmann and McKay~\cite{plantri3,plantri1,plantri2}.
Additionally, we ask how to compute construction sequences for a given planar graph, and for an improved understanding of the influence of the insertion degree.
We remark that, assuming we already have a cycle-inducing tree decomposition, our proof yields an $ O(k\cdot n) $-time algorithm for computing a construction sequence for an $n$-vertex triangulation, which we discuss more detailed in \cref{se:open_questions}.

\subsection{Preliminaries}\label{se:preliminaries}

All our graphs are finite, without loops and multi-edges (except for planar duals).
For a graph $G$ we denote its vertex-set as $V(G)$ and its edge-set as $E(G)$.

Let $H$ be a planar graph with a fixed planar embedding, called a \emph{plane} graph, and $G \subseteq H$ be a subgraph of $H$ together with its inherited planar embedding.
Then, for a face $f$ of $G$ we define $\interior_H(f)$ to be the set of all vertices in $H$ that are embedded in the strict interior of $f$. 
Note that this is the empty set when $f$ is a face of $H$, so, in particular, when $G = H$.
Furthermore, we define $\interior_H[f]$ to be the union of $\interior_H(f)$ and the set of all vertices in $H$ that are incident to $f$.
Similarly, for an induced cycle $C\subseteq H$, we define $\interior_H[C]$ to be the set of vertices in $H$ in the interior of $C$ including $V(C)$.
For a graph $G$ and a set of vertices $A \subseteq V(G)$ we denote by $G[A]$ the subgraph induced by $A$.
If $G$ is considered with a fixed planar embedding, then $G[A]$ inherits an embedding from $G$.
For a face $f$, we refer to the number of vertices that are incident to $f$ as the \emph{size} of $f$, which we denote by $ |f| $.

\subparagraph{Treewidth.} 
Let $G$ be a graph.
A tree decomposition $\td$ of $G$ is a tree whose vertices, so-called \emph{bags}, are subsets of $V(G)$ such that the following three conditions are fulfilled.
\begin{enumerate}
    \item Every vertex $v \in V(G)$ is in at least one bag.
    \item For every edge $uv \in E(G)$ there is a bag $X \in V(\td)$ such that $u,v \in X$.
    \item For every vertex $v\in V(G)$ the subgraph of $\td$ that is induced by the bags containing $v$, i.\,e., the graph $\td[\{X \in V(\td) \colon v \in X\}]$, is connected.
\end{enumerate}
The width of $\td$ is $\max\{|X| - 1 \colon X \in V(\td)\}$.
The \emph{treewidth} of $G$ is the minimum width over all tree decompositions of $G$.
We remark that for ease of presentation, we identify the vertices of a tree decomposition with the corresponding bags. Technically, this might result in the bag-set being a multi-set. However, we make sure that, except for some intermediate steps while modifying a tree decomposition, this does not occur. Hence, if two bags are equal, they refer to the same vertex of the tree.
Let $\td$ be a tree decomposition of a graph $G$. We introduce the following notation:
\begin{itemize}
    \item Let $\td' \subseteq \td$. Then, we define $G[\td']$ to be the graph induced by the set $\bigcup_{X\in V(\td')} X$.
    \item For an edge $XY \in E(\td)$ we refer to $X \cap Y$ as the \emph{adhesion set} of $XY$.
\end{itemize}

\section[Planar Quasi-\texorpdfstring{$k$}{k}-Trees]{Planar Quasi-\texorpdfstring{$\bm{k}$}{k}-Trees} \label{sec:quasi_trees}

In this section, we define planar quasi-$k$-trees and prove \cref{thm:quasi_k_trees}.
The first part of the proof is to show that planar quasi-$k$-trees are planar graphs of treewidth at most $ k $.
The main challenge is then the second part for which we show that every planar graph of treewidth at most $k$ is a subgraph of a planar quasi-$k$-tree.
Our definition generalizes the well-known construction of planar $3$-trees.
Intuitively speaking, we start with a cycle of length at most $k$ and then repeatedly choose a face, insert a vertex into it and connect it to at least one vertex that is incident to the face while ensuring that all faces have size at most $k$.
So the differences to the construction of planar $3$-trees are, first, that faces may be larger than triangles, and second, that we are allowed to omit some edges when we add a new vertex.
Note that the latter only makes a difference for $ k \geq 4 $ since for $ k = 3 $, to ensure that every face has size at most $3$, we need to add all three possible edges.
Thus, planar quasi-$3$-trees are exactly the planar $3$-trees.
Formally, we define planar quasi-$k$-trees as follows.
For this, for a graph $G$ and a vertex-ordering $\sigma= (v_1,\dots,v_n)$ of $G$, we refer to $G[\{v_1,\dots,v_i\}]$ as $G_i^\sigma$.

\begin{definition}\label{def:quasi_k-trees}
    Let $k \geq 3$ and let $ G $ be a planar graph with the following properties:
    There is a planar embedding, a vertex-ordering $\sigma = (v_1,\dots,v_n)$ and an integer $\ell\leq k$ such that
    \begin{itemize}
        \item The graph $G^\sigma_\ell$ is a cycle that bounds the outer face of $G$,
        \item for every $i\leq n$, each face in $G^\sigma_i$ has size at most $k$, and
        \item for every $i\leq n$, it holds that $G^\sigma_i$ is connected.
    \end{itemize}
    Then, we call $G$ a \emph{planar quasi-$k$-tree} and $\sigma$ a \emph{$k$-construction sequence} of $G$.
    We call a planar embedding of $G$ for which a $k$-construction sequence exists a \emph{canonical} embedding.
\end{definition}

We remark that instead of the third property, we may also require that $\deg_{G^\sigma_i}(v_i) \geq 1$ for every $1 < i \leq n$, i.\,e., whenever we insert a vertex we add at least one edge incident to it.
Furthermore, we remark that for a biconnected planar quasi-$k$-tree instead of adding singular vertices, we could in each step either add a vertex with degree at least~$2$ or choose a face bounded by vertices $(u_1,\dots,u_\ell)$ and add an induced path $(p_1,\dots,p_{\ell'})$ with $\ell + \ell' \leq k+1$ that is adjacent to exactly two vertices $u_i \neq u_j$ with $|i-j| \neq 1$.
That is, at each step of the construction sequence, at most one face is bounded by a cycle with a single induced path attached to it and each other face is bounded by a cycle.

Note that planar quasi-$k$-trees are planar by definition, so to prove the first part of \cref{thm:quasi_k_trees}, it is left to bound their treewidth by $k$.

\begin{lemma}\label{le:treewidth}
    Every planar quasi-$k$-tree $G$ has treewidth at most $k$.
\end{lemma}

\begin{proof}
    Let $\sigma=(v_1,\dots,v_n)$ be a $k$-construction sequence of $G$.
    We construct a tree decomposition maintaining that in each step of the construction, for each face there is a bag containing all vertices on the boundary.
    For the base case, take a bag with the initial cycle bounding the outer face of $ G $.
    Then, for each vertex $ v_i $ that is inserted into a face $ f_i $, append a bag containing $ v_i $ and the vertices of $ f_i $ to the already existing bag with $ f_i $.
    Since each face that occurs during the construction has size at most $ k $, the bags have size at most $ k + 1 $.
\end{proof}

We say that a tree decomposition $\td$ of a graph $G$ is \emph{cycle-inducing} if the adhesion set of each edge $XY \in E(\td)$ induces a cycle in $G$, i.\,e., $G[X \cap Y]$ is a cycle. 
Next, we use these tree decompositions to show that every planar triangulation of treewidth $k$ is a planar quasi-$k$-tree.
To do this, we need the following theorem.

\cycleinducing*

As the proof of \cref{thm:nice_edges} is fairly technical, we first demonstrate in this section how to use it, and postpone its proof until \cref{se:cycle_inducing_tree_decompositions}.
That is, for now we assume \cref{thm:nice_edges}  holds.
To start with, we make the following easy observation.

\begin{observation}\label{obs:induced_tree_decompositions}
    Let $\td$ be a cycle-inducing tree decomposition of a graph $ G $. Then, each subtree $\td' \subseteq \td$ is a cycle-inducing tree decompositions of the graph $G[\td']$.
\end{observation}

As our goal is to do induction on a cycle-inducing tree decomposition of a triangulation, we want to direct the underlying tree in some sense.
To do this, given a plane triangulation $G$ with a cycle-inducing tree decomposition $\td$, we look at bags that contain the outer triangle of $G$.
First, there is at least one such bag as each clique occurs in a bag.
We argue that if there are multiple bags with the outer triangle, one of them is a leaf in $\td$.
So consider two adjacent bags $X$ and $Y$ containing the outer triangle.
Since also $X \cap Y$ contains the outer triangle and induces a cycle in $G$ by the definition of cycle-inducing tree decompositions, $X\cap Y$ is exactly the set of vertices incident to the outer triangle.
Moreover, as the outer triangle is not a separator, one of the components of $\td - XY$ contains only the vertices in $X\cap Y$.
Thus, if there are multiple bags that contain the outer triangle, there is a leaf-bag that contains exactly the outer triangle.

Given this, we define the \emph{root} of $\td$ either as the unique bag that contains the outer triangle of $G$ or some leaf-bag containing exactly the outer triangle of $G$.
Subtrees $\td' \subseteq \td$ inherit the root, i.\,e., the root of $\td'$ is defined as the bag of $\td'$ that is in $\td$ closest to the root of $\td$. 
For a non-root bag $X\in V(\td)$ we refer to the neighbor of $X$ that is closest to the root of $\td$ as the \emph{parent} of $X$.
To all other neighbors of $X$ we refer to as the \emph{children} of $X$.

Before proving the remaining direction of \cref{thm:quasi_k_trees}, we prepare two lemmas.
First, we consider the components obtained by deleting an edge from a cycle-inducing tree decomposition.
Then, we show that the root of a cycle-inducing tree decomposition of a triangulation behaves nicely.
Both lemmas are essential for the upcoming induction in \cref{le:quasi_k-trees}.

\begin{lemma}\label{cor:inner_outer}
    Let $G$ be a plane triangulation with a cycle-inducing tree decomposition $\td$, let $XY \in E(\td)$ be an edge of $\td$ and $C$ the cycle induced by $ X \cap Y $.
    Then, the component of $\td - XY$ containing the root of $\td$ is a cycle-inducing tree decomposition of $C$ and its exterior, and the other component is a cycle-inducing tree decomposition of $C$ and its interior.
\end{lemma}

\begin{proof}
    Let $\td_1$ and $\td_2$ be the components of $\td-XY$.
    By \cref{obs:induced_tree_decompositions}, $\td_1$ and $\td_2$ are cycle-inducing tree decompositions.
    Observe that any two vertices of $ G $ that are on the same side of $C$ are not separated by $ C $ since $ G $ is a triangulation.
    As $G[\td_1]-C$ and $G[\td_2]-C$ are separated by $C$, the interior, respectively exterior, of $ C $ contains only vertices of $G[\td_1]$ or $G[\td_2]$, respectively.
    Since the root of $\td$ contains the vertices incident to the outer face of $G$, we have that the component that contains the root is a tree decomposition of $C$ and its exterior, while the other one is a tree decomposition of $C$ and its interior.
\end{proof}
\begin{restatable}{lemma}{cycleproperties} \label{le:outer_face_in_root} \label{le:small_faces}
    Let $G$ be a plane triangulation with a cycle-inducing tree decomposition $\td$ of width $k \geq 3$ and $\td' \subseteq \td$ be a subtree with root $X \in V(\td')$.
    Then, the following holds:
    \begin{enumerate}[(a)]
        \item The plane graph $G[\td']$ is bounded by an induced cycle $C$ with $C \subseteq G[X]$ and
        \item each face of $G[X]$ is bounded by a cycle of length at most $k$.
    \end{enumerate}
\end{restatable}

\begin{proof}
    First, we show that $G[\td']$ is bounded by an induced cycle contained in $G[X]$.
    If $X$ is the root of $\td$, it contains the triangle bounding $G$ (and thus also bounding $G[\td']$) by definition of the root.
    Otherwise, let $Y\in V(\td)$ be the parent of $X$ in $\td$ and consider the cycle $C = G[X \cap Y]$.
    By \cref{cor:inner_outer}, $G[\td']$ is bounded by $C$ as $\td'$ does not contain the root.
    Clearly, it holds that $C \subseteq G[X] $.
    
    It is left to show that each face of $G[X]$ is bounded by a cycle of length at most $k$.
    For this we assume without loss of generality that there are no adjacent bags in $\td$ that contain the same set of vertices.
    If there are two such bags, we can contract the edge connecting them and argue on the resulting tree decomposition.
    Recall that the adhesion set of an edge $AB \in E(\td)$ is defined as $A \cap B$.
    Since we assume that $A \neq B$ for every $AB \in E(\td)$, this gives us that the adhesion set of each edge in $\td$ has size at most $k$.
    That is, to prove that each face of $G[X]$ is bounded by a cycle of length at most $k$, it suffices to show that each face of $G[X]$ is either a triangle or bounded by a cycle induced by the adhesion set of an edge of $\td$.
    
    For this let $f$ be a face of $G[X]$.
    If $f$ is a triangle we are done.
    Hence, assume in the following that $f$ is incident to at least four vertices.
    We aim to find a neighbor $ Z $ of the root $ X $ such that $ XZ $ induces the boundary of $ f $.
    As $G$ is a triangulation, there is a vertex $v\in \interior_G(f)$, i.\,e., in the interior of $f$ in $G$ and in particular $ v \notin X $.
    Consider the component $\td'' \subseteq \td - X$ such that $v$ is in some bag of $\td''$ and let $Z$ be the neighbor of $X$ that is in $\td''$. 
    Observe that since $v\notin X$, there is indeed only one such component.
    Now consider the edge $XZ$ and the cycle $ C = G[X \cap Z] $ it induces.
    It is left to show that $ C $ indeed bounds $ f $.
    First, the vertices $ V(f) $ on the boundary of $ f $ are in $ \interior_G[C] $ since the interiors of $ f $ and $ C $ share a vertex, namely $ v $, and as a face of $ G[X] $, $ f $ does not contain any vertices of $ C \subseteq G[X] $ in its strict interior.
    Second, no vertex of $ V(f) $ is in the strict interior of $ C $ since by \cref{cor:inner_outer}, no vertex of $ X $ is in the strict interior of $ C $, so, in particular, no vertex of the face $ f $ of $ G[X] $.
    Together, we have that the boundary of $ f $ is a subset of $ C $, and since $ C $ is an induced cycle, we have that $ C $ contains all vertices of $ V(f) $.
    To sum up, $ f $ is bounded by the cycle $ C $, which is induced by an adhesion set and, thus, has size at most $ k $.
\end{proof}

By \cref{thm:nice_edges}, every triangulation of treewidth $k$ admits a cycle-inducing tree decomposition of the same width. Next, we use this to show that every triangulation of treewidth $k$ is a planar quasi-$k$-tree.
To do so, we show the following stronger statement by induction.

\begin{lemma} \label{le:quasi_k-trees}
    Let $G$ be a plane triangulation with a cycle-inducing tree decomposition $\td$ of width $k$. Then, for all $\td' \subseteq \td$ the graph $G[\td']$ is a planar quasi-$k$-tree and the embedding of $G[\td']$ is canonical.
\end{lemma}

\begin{proof}
    We fix $ \td $ and show the statement by induction on the number $|V(\td')|$ of bags of $ \td' $, so we aim for a $ k $-construction sequence of $ G[\td'] $.
    We start with some observations that are useful both for the base case and the step.
    For this, let $X \in V(\td')$ be the root of $\td'$.
    By \cref{le:small_faces}a, we know that $G[X]$ is bounded by an induced cycle of length at most $k$.
    Furthermore, by \cref{le:small_faces}b, we have that $G[X]$ is connected as every face of $G[X]$ is bounded by a cycle.
    Let $\sigma = (v_1,\dots,v_{|X|})$ be some ordering of $X$ such that $G_i^\sigma$ is connected for all $i \leq |X|$ and $v_1,\dots,v_\ell$ are the vertices incident to the outer face of $G[X]$ for some $\ell \leq k$.
    Then, for all $i \leq k$, each face of $G^\sigma_i$ contains at most $i\leq k$ vertices and by the definition of $\sigma$, every $G^\sigma_i$ is connected.
    Furthermore, placing $v_{|X|}$ does not create a face of size $k+1$, as $G^\sigma_{|X|}=G[X]$ does not contain such a face by \cref{le:small_faces}b.
    Thus, $\sigma$ is a $k$-construction sequence of $G[X]$, $G[X]$ is a planar quasi-$k$-tree, and the fixed embedding of $G[X]$ is a canonical embedding.
    Therefore, if $|V(\td')| = 1$ we are done.
    
    Next, assume that $|V(\td')| > 1$,
    let $d = \deg_{\td'}(X)$ and $Y_1,\dots, Y_d$ be the children of $X$ in $\td'$.
    By \cref{cor:inner_outer}, the components $\td_1,\dots,\td_d\subseteq \td' - X$ are cycle-inducing tree decompositions of the cycle $ G[X \cap Y_i] $ together with its interior.
    Thus, by induction there are $k$-construction sequences $\sigma_1,\dots \sigma_d$ of $G[\td_1], \dots, G[\td_d]$, respectively.
    Furthermore, the embedding of $G$ restricted to $G[\td_i]$ is canonical for all $i \leq d$.

    \begin{figure}
        \centering
        \includegraphics{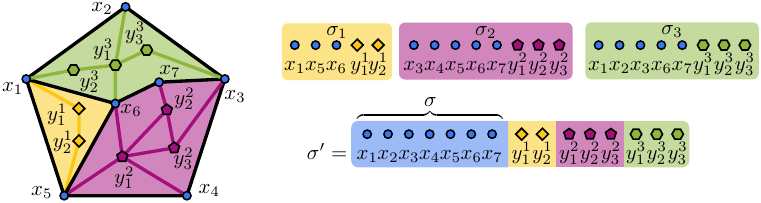}
        \caption{Example for the construction of $\sigma'$
        with $\sigma = (x_1,\dots,x_7)$ being a $6$-construction sequence of $G[X]$ and $\sigma_1,\sigma_2$ and $\sigma_3$ being $6$-construction sequences of $G[\td_1],G[\td_2]$ and $G[\td_3]$, respectively.}
        \label{fig:induction_step}
    \end{figure}

    To show that $G[\td']$ is a planar quasi-$k$-tree, it is left to combine $\sigma_1,\dots \sigma_d$ together with the $k$-construction sequence $\sigma$ of $G[X]$ to a $k$-construction sequence of $G[\td']$; for an example see \cref{fig:induction_step}.
    For this, recall that the vertices in $ X \cap Y_i $ are the outer vertices of $ G[\td_i] $, while all other vertices of $ G[\td_i] $ are in the strict interior of the cycle $ G[X \cap Y_i] $.
    By the definition of $ k $-construction sequences, this implies that $\sigma_i$ starts with the vertices in $X \cap Y_i$, and that the $ V(G[\td_i]) $ are pairwise disjoint except for vertices in the root $ X $.
    Now let $\sigma_1',\dots,\sigma_d'$ be $\sigma_1,\dots,\sigma_d$ restricted to the vertices $V(G[\td']) - X$, i.\,e., without the outer vertices at the beginning of each sequence.
    Then, $ \sigma' = \sigma\cdot \sigma_1'\cdot \ldots \cdot \sigma_d'$ contains each vertex of $ G[\td'] $ exactly once and is a $k$-construction sequence:
    The first property, i.\,e., that the sequence starts with the cycle bounding the outer face of $ G[\td'] $ is met due to the definition of $ \sigma $, and the second and third property are inherited from the subsequences.
    Hence, the fixed embedding of $ G[\td'] $ is canonical.
\end{proof}

With \cref{le:quasi_k-trees} we have that triangulations of treewidth at most $k$ are planar quasi-$k$-trees. 
To show that every planar graph of treewidth at most $k$ is a subgraph of a planar quasi-$k$-tree, we use the following well-known result by Bield and Ruiz Velázquez.

\begin{lemma}[Biedl and Ruiz Velázquez~\cite{triangulations_with_same_tw}] \label{le:triangulations_with_same_tw}
    For every planar graph $G$ of treewidth~$k$, there is a planar triangulation $H$ of $G$ on the same vertex-set such that $H$ has treewidth $\max\{3,k\}$.
\end{lemma}

\begin{corollary}\label{cor:quasi_k-trees_vertex-set}
    For every planar graph $G$ of treewidth $k \geq 3$, there exists a planar quasi-$k$-tree $H$ with $G \subseteq H$ on the same vertex-set.
\end{corollary}

\begin{proof}
    Let $G$ be a planar graph of treewidth at most $k$.
    By \cref{le:triangulations_with_same_tw}, there exists a triangulation $H$ of $G$ of treewidth at most $k$ and the same set of vertices as $G$. 
    Thus, by \cref{le:tree_decomposition_nice_edges}, there exists a cycle-inducing tree decomposition of $H$ of width $k$, giving us that $H$ is a planar quasi-$k$-tree by \cref{le:quasi_k-trees}.
\end{proof}
\cref{le:treewidth,cor:quasi_k-trees_vertex-set} together prove both directions of \cref{thm:quasi_k_trees}, assuming \cref{thm:nice_edges}.

\section{Planar Quasi-4-Trees and Vertex Arboricity}\label{se:quasi_4-trees}

In this section, we prove \cref{thm:arboricity,thm:nice_4-trees}.
\cref{thm:arboricity} states that the vertex-set of every planar graph $G$ of treewidth at most $4$ can be partitioned into two sets, each of which induces a forest. 
In other words, the vertex arboricity of $G$ is at most~$2$.
To show this, we first define \emph{\nice} planar quasi-$k$-trees and show that \cref{thm:nice_4-trees} holds.
That is, we show that every planar graph of treewidth at most $4$ is a subgraph of a \nice planar quasi-$k$-tree.

We say that a planar quasi-$k$-tree is a \emph{\nice} planar quasi-$k$-trees if the graph is biconnected at every step of the construction.
For $k=4$ this yields exactly planar quasi-$4$-trees as originally defined by Förster~\cite{planar_quasi_4_trees}.
Förster conjectured that every planar graph of treewidth at most $4$ is a subgraph of a \nice planar quasi-$4$-tree.
We confirm this by proving \cref{thm:nice_4-trees}.
Formally, we define \nice planar quasi-$k$-trees as follows.

\begin{definition}\label{def:strong_quasi_k-trees}
    Let $G$ be a planar quasi-$k$-tree.
    If there is a $k$-construction sequence $\sigma= (v_1,\dots,v_n)$ of $G$ such that $G^\sigma_i$ is biconnected for all $3\leq i \leq n$, we call $\sigma$ a \emph{\nice} $k$-construction sequence of $G$ and $G$ a \emph{\nice} planar quasi-$k$-tree.
\end{definition}

Let us remark that in a recent book chapter on product structure Dujmovi\'{c} introduced \emph{plane $(k,s)$-trees}~\cite{plane_ks_trees}.
For $s=1$ these yield almost the same graph class as \nice planar quasi-$k$-trees.

For an alternative view of \cref{def:strong_quasi_k-trees}, given a graph $G$ and a $k$-construction sequence $\sigma=(v_1,\dots,v_n)$ of $G$, we define the \emph{insertion degree} of a vertex $v_i$ as $\deg_{G^\sigma_i}(v_i)$, i.\,e., the degree of $v_i$ at the time of its insertion.
Note that a $k$-construction sequence of a graph $G$ is a \nice $k$-construction sequence if and only if every vertex (except for the first two) has insertion degree at least $2$.

Since triangulations have minimum degree larger than 1, the following lemma together with \cref{thm:quasi_k_trees} implies \cref{thm:nice_4-trees}, which states that for planar graphs, treewidth at most 4 is equivalent to being a subgraph of a strong planar quasi-4-tree.
We remark that the requirement of minimum degree at least 2 is necessary as otherwise the graph is not biconnected.
To complement \cref{thm:nice_4-trees}, we show in \cref{sec:niceness} that the same does not hold for $k > 4$ by constructing a graph with treewidth $5$ that is not a subgraph of any \nice planar quasi-$k$-tree.

\begin{lemma}
    Let $G$ be a planar quasi-$4$-tree with minimum degree $\delta(G) > 1$. Then, $G$ is a \nice planar quasi-$4$-tree.
\end{lemma}

\begin{figure}
    \centering
    \includegraphics{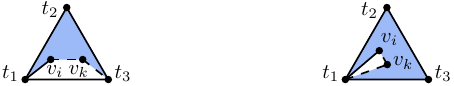}
    \caption{The two options (up to isomorphism) for the graph induced by $t_1,t_2,t_3,v_j,$ and $v_k$. Both contain a face of size $5$ (blue).
    One of the dashed edges in each graph might not exist.}
    \label{fig:face_size_5}
\end{figure}

\begin{proof}
    Among all 4-construction sequences of $ G $, let $\sigma=(v_1,\dots,v_n)$ be one such that the first vertex with insertion degree 1 occurs as late as possible, if at all. 
    Assume there is such a vertex $ v_i $ in $\sigma$, i.\,e., $ v_i $ has insertion degree 1 but the insertion degree of all $ v_{j} $ with $ j < i $ is at least 2.
    Note that $v_i$ is inserted into a triangle in $G^\sigma_{i-1}$ as otherwise the insertion of $ v_i $ creates a face of size larger than 4.
    Let $t_1,t_2,$ and $t_3$ be the vertices of this triangle. 
    
    Now, as $\delta(G)>1$, there is a vertex $v_k$ with $k > i$ that is inserted into the face bounded by $t_1,t_2,t_3,$ and $v_i$. 
    Observe that $v_k$ is adjacent to at most one of $t_1,t_2,$ and $t_3$, as otherwise swapping $v_i$ and $v_k$ in $ \sigma $ yields a 4-construction sequence that is biconnected until the insertion of $ v_k $, contradicting the choice of $ \sigma $. 
    Hence, the graph $G^\sigma_k$ contains a face of size $5$ (see \cref{fig:face_size_5}). 
    Therefore, there is no such vertex $v_i$ and $G$ is indeed a \nice planar quasi-$4$-tree.
\end{proof}

Now, to show that planar graphs of treewidth 4 have vertex arboricity at most~$2$, we may assume that $G$ is a strong planar quasi-$4$-tree by \cref{thm:nice_4-trees}, since the vertex arboricity is monotone under taking subgraphs.
Even stronger, by additionally using \cref{le:triangulations_with_same_tw}, we may assume that $G$ is a triangulated strong planar quasi-$4$-tree. 
By a theorem of Stein~\cite{characterization_arboricity}, a planar triangulation $G$ has vertex arboricity~$2$ if and only if the dual of $G$ is Hamiltonian.

\begin{theorem}[Stein, 1971~\cite{characterization_arboricity}]\label{thm:stein}
    Let $G$ be a planar triangulation with $|V(G)| \geq 4$.
    Then, $V(G)$ can be partitioned into two induced forests if and only if the dual of $G$ admits a Hamiltonian cycle.
\end{theorem}

Hence, given a triangulated \nice planar quasi-$4$-tree $G$, our task is to find a Hamiltonian cycle $C$ in the dual $G^*$ of $G$.
To do so, we proceed along a strong $4$-construction sequence $\sigma = (v_1,\ldots,v_n)$ corresponding to a canonical embedding of $G$ (cf.\ \cref{def:quasi_k-trees,def:strong_quasi_k-trees}) and aim for a closed walk visiting all vertices.
That is, in contrast to a Hamilton cycle, we may see vertices multiple times.
For an edge $e \in E(G)$ we denote by $e^*$ the dual edge of $e$ in $G^*$. 

\begin{figure}
        \centering
        \begin{subfigure}[t]{0.3\textwidth}
            \centering
            \includegraphics[page=1]{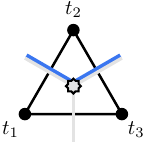}
            \subcaption{}
            \label{fig:interaction_triangle}
        \end{subfigure}%
        \hfill
        \begin{subfigure}[t]{0.3\textwidth}
            \centering
            \includegraphics[page=2]{figures/hamilton_cycle_subfigures.pdf}
            \subcaption{}
            \label{fig:interaction_quadrangle_two_edges}
        \end{subfigure}%
        \hfill
        \begin{subfigure}[t]{0.3\textwidth}
            \centering
            \includegraphics[page=3]{figures/hamilton_cycle_subfigures.pdf}
             \subcaption{}
             \label{fig:interaction_quadrangle_four_edges}
        \end{subfigure}%
        \caption{
            The three possible ways for the closed walk $C$ (blue and orange) to interact with a face of $G$ (black).
            The dual of $G$ is depicted in gray.}
        \label{fig:allowed_patterns}
    \end{figure}

\begin{restatable}{lemma}{arboricity}\label{le:arboricity}
    Let $G$ be a \nice planar quasi-$4$-tree with a canonical embedding and $G^*$ be the dual of $G$.
    Then, there is a closed walk $C$ in $G^*$ that fulfills the following properties:
    \begin{itemize}
        \item Every edge in $E(G^*)$ is contained at most once in $C$.
        \item Every vertex in $V(G^*)$ is contained at least once in $C$.
        \item Each pair of consecutive edges $uv$ and $vw$ in $C$ is consecutive in the cyclic ordering of edges incident to $v$.
    \end{itemize}
    
\end{restatable}
Before we prove \cref{le:arboricity}, let us observe some properties of such a closed walk $C$ in $G^*$.
\cref{fig:allowed_patterns} shows the possible interactions of $C$ with a facial cycle of $G$.
Recall for this that for every $i \in \{3,\ldots,n\}$, the graph $G_i^\sigma = G[\{v_1,\ldots,v_i\}]$ is biconnected, and each face of $G_i^\sigma$ is incident to $3$ or $4$ vertices.
In particular, each triangular face of $G$ contains exactly two dual edges of $C$, while a quadrangular face of $G$ contains either two or four dual edges of $C$.
Note that for a triangulation $G$, \cref{le:arboricity} yields a Hamiltonian cycle in $G^*$.
That is, by \cref{thm:stein} and \cref{le:triangulations_with_same_tw}, \cref{le:arboricity} implies \cref{thm:arboricity}.

\begin{proof}
    We proceed by induction on the number $n$ of vertices in $G$.
    For $n \leq 3$, $C$ is formed by taking any two edges of $G^*$.
    Thus, let $n > 3$ and consider the last vertex $v$ in a \nice $4$-construction sequence of $G$.
    Then, $G-v$ is a \nice quasi-$4$-tree inheriting a canonical embedding from the canonical embedding of $G$.
    Let $f$ be the face of $G-v$ that contains $v$ in its interior.
    By the definition of \nice construction sequences, $v$ has degree at least $2$ in $G$.
    
    In the following, we do a case distinction based on whether $f$ is triangle or quadrangle, and whether $v$ is in $G$ adjacent to two, three, or four vertices of $f$.
    In almost all cases we apply induction to the \nice quasi-$4$-tree $G' = G-v$, obtain a closed walk $C'$ in the dual of $G'$ (as depicted in \cref{fig:allowed_patterns}), and do local modifications (shown in \cref{fig:triangle_cases,fig:four_edges,fig:two_edges}) to obtain the desired closed walk $C$ in the dual of $G$.
    Only in case $f$ is a quadrangle and $v$ has three neighbors on $f$ (which we handle at the end of the proof), we add one more edge to $G'$ before calling induction.

    \medskip

    \begin{figure}
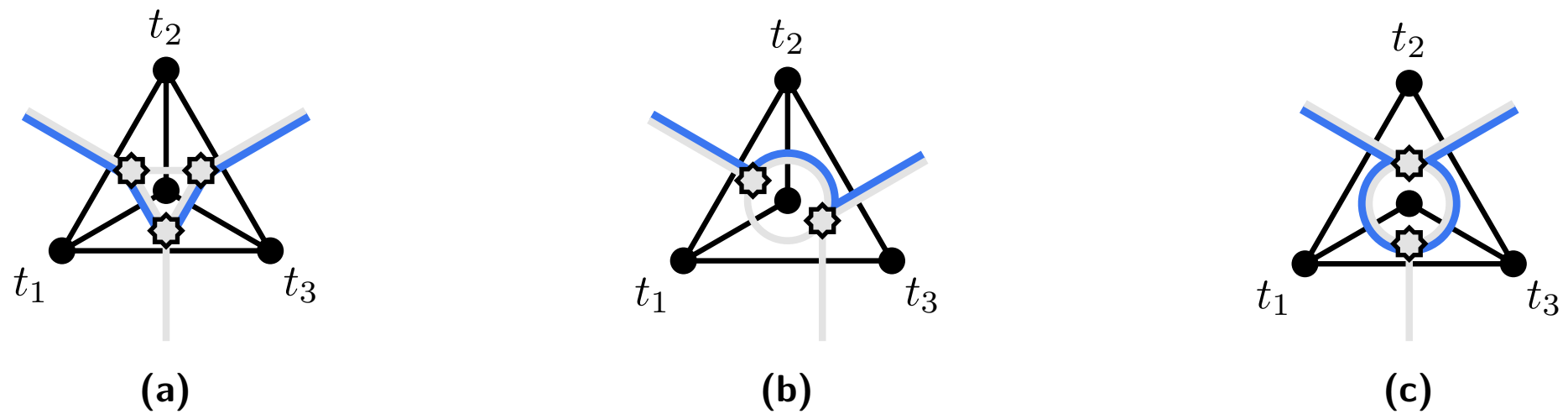

        \centering
        \begin{subfigure}[t]{0.3\textwidth}
            \centering
            \includegraphics[page=4]{figures/hamilton_cycle_subfigures.pdf}
            \subcaption{}
            \label{fig:triangle_case_123}
        \end{subfigure}%
        \hfill
        \begin{subfigure}[t]{0.3\textwidth}
            \centering
            \includegraphics[page=5]{figures/hamilton_cycle_subfigures.pdf}
            \subcaption{}
            \label{fig:triangle_case_12}
        \end{subfigure}%
        \hfill
        \begin{subfigure}[t]{0.3\textwidth}
            \centering
            \includegraphics[page=6]{figures/hamilton_cycle_subfigures.pdf}
             \subcaption{}
             \label{fig:triangle_case_13}
        \end{subfigure}%
        \caption{
            The three cases if $f$ is bounded by a triangle.
            The augmented closed walk $C$ is highlighted in blue.
            The dual of $G$ is depicted in gray.}
        \label{fig:triangle_cases}
    \end{figure}

    \medskip

    \textbf{Case 1.} \ 
    Suppose $f$ is a triangle $\{t_1,t_2,t_3\}$.
    Let $C'$ be a closed walk in the dual of $G' = G-v$ obtained by induction.
    Without loss of generality $C'$ contains the edges of $(t_1t_2)^*$ and $(t_2t_3)^*$.
    Vertex $v$ is adjacent to at least two of $t_1,t_2$ and $t_3$.
    If $v$ is adjacent to $t_1,t_2$ and $t_3$, as depicted in \cref{fig:triangle_case_123}, we take $C = C' + \{(t_1v)^*,(t_3v)^*\}$.
    If $v$ is adjacent to $t_1$ and $t_2$, as depicted in \cref{fig:triangle_case_12}, we take $C = C' + \{(t_2v)^*\}$.
    Finally, if $v$ is adjacent to $t_1$ and $t_3$, as depicted in \cref{fig:triangle_case_13}, we take $C = C' + \{(t_1v)^*,(t_3v)^*\}$.

    \textbf{Case 2.} \ 
    Suppose $f$ is a quadrangle $\{q_1,q_2,q_3,q_4\}$, and that $v$ is incident to exactly two or four vertices of $f$.
    Let $C'$ be a closed walk in the dual of $G' = G-v$ obtained by induction.
    Without loss of generality $C'$ contains the edges $(q_1q_2)^*$ and $(q_2q_3)^*$ consecutively.
    
    Suppose further, \textbf{Case 2.1.}, that $C'$ also contains the edges $(q_3q_4)^*$ and $(q_4q_1)^*$.
    If $v$ is adjacent to $q_1,q_2,q_3$ and $q_4$, as depicted in \cref{fig:quadrangle_case_four_edges1234}, we take $C = C' + \{(q_2v)^*,(q_4v)^*\}$.
    If $v$ is adjacent to $q_1$ and $q_3$, as depicted in \cref{fig:quadrangle_case_four_edges13}, we take $C = C'$.
    If $v$ is adjacent to $q_2$ and $q_4$, as depicted in \cref{fig:quadrangle_case_four_edges24}, we take $C = C' + \{(q_2v)^*,(q_4v)^*\}$.

    \begin{figure}
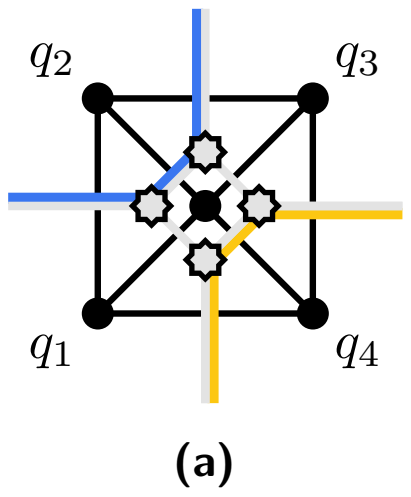

        \centering
        \begin{subfigure}[t]{0.3\textwidth}
            \centering
            \includegraphics[page=7]{figures/hamilton_cycle_subfigures.pdf}
            \subcaption{}
            \label{fig:quadrangle_case_four_edges1234}
        \end{subfigure}%
        \hfill
        \begin{subfigure}[t]{0.3\textwidth}
            \centering
            \includegraphics[page=8]{figures/hamilton_cycle_subfigures.pdf}
            \subcaption{}
            \label{fig:quadrangle_case_four_edges13}
        \end{subfigure}%
        \hfill
        \begin{subfigure}[t]{0.3\textwidth}
            \centering
            \includegraphics[page=9]{figures/hamilton_cycle_subfigures.pdf}
             \subcaption{}
             \label{fig:quadrangle_case_four_edges24}
        \end{subfigure}%
        \caption{
            The three cases if $f$ is bounded by a quadrangle and the closed walk $C'$ contains all dual edges of edges incident to $f$.
            The augmented closed walk $C$ is highlighted in blue and orange.
            The dual of $G$ is depicted in gray.}
        \label{fig:four_edges}
    \end{figure}

    Now suppose, \textbf{Case 2.2.}, that $C'$ does not contain the edges $(q_3q_4)^*$ and $(q_4q_1)^*$.
    If $v$ is adjacent to $q_1,q_2,q_3$ and $q_4$, as depicted in \cref{fig:quadrangle_case_two_edges1234}, we take $C = C' + \{(q_1v)^*,(q_4v)^*,(q_3v)^*\}$.
    If $v$ is adjacent to $q_1$ and $q_3$, as depicted in \cref{fig:quadrangle_case_two_edges13}, we take $C = C' + \{(q_1v)^*,(q_3v)^*\}$.
    If $v$ is adjacent to $q_2$ and $q_4$, as depicted in \cref{fig:quadrangle_case_two_edges24}, we take $C = C' + \{(q_2v)^*\}$.

    \begin{figure}
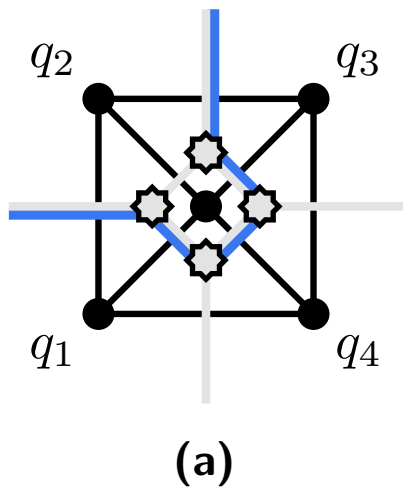

        \centering
        \begin{subfigure}[t]{0.3\textwidth}
            \centering
            \includegraphics[page=10]{figures/hamilton_cycle_subfigures.pdf}
            \subcaption{}
            \label{fig:quadrangle_case_two_edges1234}
        \end{subfigure}%
        \hfill
        \begin{subfigure}[t]{0.3\textwidth}
            \centering
            \includegraphics[page=11]{figures/hamilton_cycle_subfigures.pdf}
            \subcaption{}
            \label{fig:quadrangle_case_two_edges13}
        \end{subfigure}%
        \hfill
        \begin{subfigure}[t]{0.3\textwidth}
            \centering
            \includegraphics[page=12]{figures/hamilton_cycle_subfigures.pdf}
             \subcaption{}
             \label{fig:quadrangle_case_two_edges24}
        \end{subfigure}%
        \caption{
            The three cases if $f$ is bounded by a quadrangle and the closed walk $C'$ contains not all dual edges of edges incident to $f$.
            The augmented closed walk $C$ is highlighted in blue.
            The dual of $G$ is depicted in gray.}
        \label{fig:two_edges}
    \end{figure}

    Suppose $f$ is a quadrangle $\{q_1,q_2,q_3,q_4\}$, and that $v$ is adjacent to exactly three vertices of $f$, say $q_1$, $q_2$, and $q_3$.

    \textbf{Case 3.} \ 
    First, we assume that $q_2q_4$ is not an edge of~$G$.
    In this case, we add the edge $q_2q_4$ before calling induction.
    Note that adding an edge to a canonical embedding of a \nice planar quasi-$4$-tree that maintains planarity results again in a \nice planar quasi-$4$-tree, as adding the edge at the time of inserting the second of the incident vertices may only decrease the size of the faces.
    That is, let $G' = G - v + \{q_2q_4\}$ and let $C'$ be a closed walk in the dual of $G'$ obtained by induction.
    In particular, the quadrangular face $f$ is split into two triangular faces $f_1 = \{q_1,q_2,q_4\}$ and $f_2 = \{q_2,q_3,q_4\}$, each of which is visited exactly once by $C'$.

    \begin{figure}
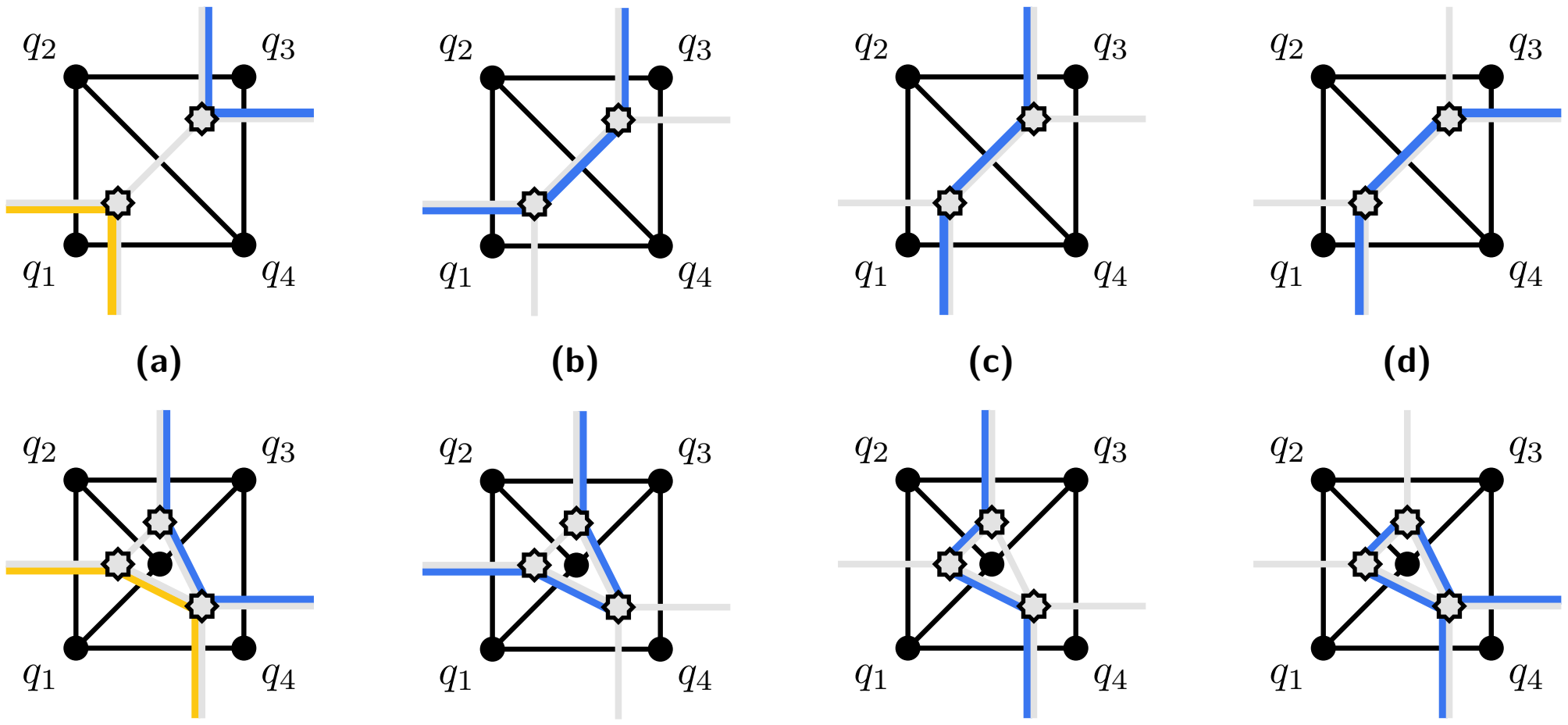

        \centering
        \begin{subfigure}[t]{0.245\textwidth}
            \centering
            \includegraphics[page=17]{figures/hamilton_cycle_subfigures.pdf}
            \subcaption{}
            \label{fig:quadrangle_case_induction_all_edges_before}
        \end{subfigure}%
        \hfill
        \begin{subfigure}[t]{0.245\textwidth}
            \centering
            \includegraphics[page=18]{figures/hamilton_cycle_subfigures.pdf}
            \subcaption{}
            \label{fig:quadrangle_case_induction1223_before}
        \end{subfigure}%
        \hfill
        \begin{subfigure}[t]{0.245\textwidth}
            \centering
            \includegraphics[page=19]{figures/hamilton_cycle_subfigures.pdf}
             \subcaption{}
             \label{fig:quadrangle_case_induction2314_before}
        \end{subfigure}%
        \hfill
        \begin{subfigure}[t]{0.245\textwidth}
            \centering
            \includegraphics[page=20]{figures/hamilton_cycle_subfigures.pdf}
             \subcaption{}
             \label{fig:quadrangle_case_induction1434_before}
        \end{subfigure}%
        \hfill

        \vspace{1ex}
        \begin{subfigure}[t]{0.245\textwidth}
            \centering
            \includegraphics[page=13]{figures/hamilton_cycle_subfigures.pdf}
            \subcaption{}
            \label{fig:quadrangle_case_induction_all_edges}
        \end{subfigure}%
        \hfill
        \begin{subfigure}[t]{0.245\textwidth}
            \centering
            \includegraphics[page=14]{figures/hamilton_cycle_subfigures.pdf}
            \subcaption{}
            \label{fig:quadrangle_case_induction1223}
        \end{subfigure}%
        \hfill
        \begin{subfigure}[t]{0.245\textwidth}
            \centering
            \includegraphics[page=15]{figures/hamilton_cycle_subfigures.pdf}
             \subcaption{}
             \label{fig:quadrangle_case_induction2314}
        \end{subfigure}%
        \hfill
        \begin{subfigure}[t]{0.245\textwidth}
            \centering
            \includegraphics[page=16]{figures/hamilton_cycle_subfigures.pdf}
             \subcaption{}
             \label{fig:quadrangle_case_induction1434}
        \end{subfigure}%
        \hfill
        \caption{
            The four cases if $f$ is bounded by a quadrangle, $v$ is incident to three vertices and the edge $q_2q_4$ is not present in $G$.
            On the top the graph $G' = G - v + \{q_2q_4\}$ with the closed walk $C'$ obtained by induction highlighted in blue and orange.
            On the bottom the graph $G$ with the augmented closed walk $C$ highlighted in blue and orange.
            The duals of $G$ and $G'$ are depicted in gray.}
        \label{fig:hamilton_induction}
    \end{figure}

    If, \textbf{Case 3.1.}, $C'$ does not visit $f_1,f_2$ consecutively (i.\,e., $C'$ does not contain $(q_2q_4)^*$), then $C'$ contains the edges $(q_1q_2)^*,(q_1q_4)^*$ consecutively and the edges $(q_2q_3)^*,(q_3q_4)^*$ consecutively, as depicted in \cref{fig:quadrangle_case_induction_all_edges_before}.
    In this case we take $C = C' + \{(q_1v)^*,(q_3v)^*\}$ as depicted in \cref{fig:quadrangle_case_induction_all_edges}.
    
    If, \textbf{Case 3.2.}, $C'$ visits $f_1,f_2$ consecutively (i.\,e., $C'$ contains $(q_2q_4)^*$), then $C'$ contains exactly one of $\{(q_1q_2)^*,(q_1q_4)^*\}$ and exactly one of $\{(q_2q_3)^*,(q_3q_4)^*\}$.
    If $C'$ contains $(q_1q_2)^*$ and $(q_2q_3)^*$ (\cref{fig:quadrangle_case_induction1223_before}), we take $C = C' + \{(q_1v)^*,(q_3v)^*\}$ as depicted in \cref{fig:quadrangle_case_induction1223}.
    If $C'$ contains $(q_1q_4)^*$ and $(q_2q_3)^*$ (\cref{fig:quadrangle_case_induction2314_before}), we take $C = C' + \{(q_1v)^*,(q_2v)^*\}$ as depicted in \cref{fig:quadrangle_case_induction2314}.
    If $C'$ contains $(q_1q_4)^*$ and $(q_3q_4)^*$ (\cref{fig:quadrangle_case_induction1434_before}), we take $C = C' + \{(q_1v)^*,(q_2v)^*,(q_3v)^*\}$ as depicted in \cref{fig:quadrangle_case_induction1434}.

    \begin{figure}
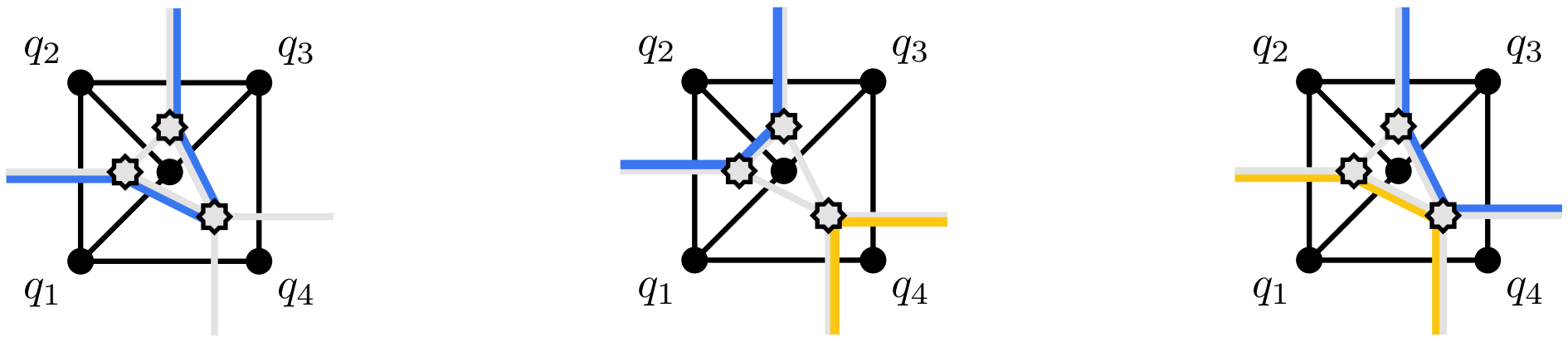

        \centering
        \begin{subfigure}[t]{0.3\textwidth}
            \centering
            \includegraphics[page=21]{figures/hamilton_cycle_subfigures.pdf}
            \subcaption{}
            \label{fig:quadrangle_case_induction1223_present}
        \end{subfigure}%
        \hfill
        \begin{subfigure}[t]{0.3\textwidth}
            \centering
            \includegraphics[page=22]{figures/hamilton_cycle_subfigures.pdf}
            \subcaption{}
            \label{fig:quadrangle_case_all_edges1_present}
        \end{subfigure}%
        \hfill
        \begin{subfigure}[t]{0.3\textwidth}
            \centering
            \includegraphics[page=23]{figures/hamilton_cycle_subfigures.pdf}
             \subcaption{}
             \label{fig:quadrangle_case_all_edges2_present}
        \end{subfigure}%
        \caption{
            The three cases if $f$ is bounded by a quadrangle, $v$ is incident to three vertices and the edge $q_2q_4$ is present in $G$.
             The augmented closed walk $C$ is highlighted in blue and orange.
            The dual of $G$ is depicted in gray.}
        \label{fig:hamilton_induction_present}
    \end{figure}

    \textbf{Case 4.} \ 
    Finally, consider the case that the edge $q_2q_4$ is an edge of~$G$.
    We assume without loss of generality that the vertices $q_1,q_2,$ and $q_4$ form a triangle with $q_3$ in its interior.
    Let $C'$ be a closed walk in the dual of $G' = G-v$ obtained by induction.
    Since there is at least one face in the interior of the triangle induced by $q_1,q_2,$ and $q_4$, the closed walk $C'$ intersects this triangle exactly twice.
    Analogously, the closed walk $C'$ intersects the triangle induced by $q_2,q_3,$ and $q_4$ exactly twice.
    That is, $C'$ contains at least one edge of $(q_2q_3)^*$ and $(q_3q_4)^*$ as well as at least one edge of $(q_1q_2)^*$ and $(q_1q_4)^*$.
    If $C'$ contains exactly $(q_1q_2)^*$ and $(q_2q_3)^*$, we take $C = C' + \{(q_1v)^*,(q_3v)^*\}$ as depicted in \cref{fig:quadrangle_case_induction1223_present}.
    If $C'$ contains $(q_1q_2)^*$ and $(q_2q_3)^*$ consecutively as well as $(q_1q_4)^*$ and $(q_3q_4)^*$ consecutively, we take $C = C' + \{(q_2v)^*\}$ as depicted in \cref{fig:quadrangle_case_all_edges1_present}.
    If $C'$ contains $(q_1q_2)^*$ and $(q_1q_4)^*$ consecutively as well as $(q_2q_3)^*$ and $(q_3q_4)^*$ consecutively, we take $C = C' + \{(q_1v)^*,(q_3v)^*\}$ as depicted in \cref{fig:quadrangle_case_all_edges2_present}.
    
    It is straightforward to verify that in every case $C$ is a closed walk in the dual of $G$ satisfying all desired properties.
\end{proof}

\section{Cycle-Inducing Tree Decompositions} \label{se:cycle_inducing_tree_decompositions}

Recall that a tree decomposition is called \emph{cycle-inducing} if the adhesion set of each edge of the tree decomposition, i.\,e., the intersection of the bags incident to the edge, induces a cycle in the graph.
In this section, we prove the following theorem, which completes the proof of \cref{thm:quasi_k_trees,thm:nice_4-trees}.

\cycleinducing*

To start with, we need the following well-known result regarding separators in triangulations.
For the sake of completeness we include a proof of this folklore result.

\begin{restatable}{lemma}{separators}\label{le:inclusion_minimal_separators_in_triangulations}
    Let $G$ be a planar triangulation and let $S\subsetneq V(G)$ be a separator of $G$. Then, $S$ induces a cycle in $G$ if and only if $S$ is an inclusion-minimal separator of $G$.
\end{restatable}

\begin{proof}
    Let $G$ be a triangulation and let $S\subsetneq V(G)$ be an inclusion-minimal separator of $G$.
    As $S$ is inclusion-minimal, every vertex $s \in S$ is adjacent to a vertex in each component of $G - S$.
    Note that $|S| \geq 3$ because $G$ is a triangulation and thus $3$-connected.
    Therefore, contracting each component of $G-S$ into a single vertex results in a $K_{3,c}$-minor with $c$ being the number of components of $G-S$.
    Since $G$ is planar, $G-S$ has exactly two components.
    
    So, let $A$ and $B$ be the two components of $G-S$ and fix some planar embedding of $G$.
    As $G[V(A) \cup S]$ is connected, it is embedded inside a single face of $B$.
    We assume without loss of generality that this face is not the outer face of $G$.
    Next, we claim that $G[S]$ contains a cycle with $A$ in its interior.
    Assume this is not the case. Then, the outer face of $G[V(A) \cup S]$ contains a vertex $a \in V(A)$.
    Thus, there is a face that contains a vertex $b \in V(B)$ and $a$ and, therefore, as $G$ is a triangulation, an edge $ab \in E(G)$.
    This contradicts $A$ and $B$ being two separate components of $G-S$.
    Thus, $S$ contains a cycle $C$ with $A$ in its interior.
    As $S$ contains such a cycle $C$, it also contains an induced cycle $C'\subseteq C$ with $A$ in its interior.
    Since $S$ is inclusion-minimal, we have that $S = V(C')$.

    For the other direction let $S\subsetneq V(G)$ be a separator of $G$ that induces a cycle.
    Since $S$ is a separator, there is an inclusion-minimal separator $S' \subseteq S$.
    Therefore, $S'$ induces a cycle in $G$.
    As $S$ induces a cycle as well we get that $S'=S$.
\end{proof}

\begin{figure}
    \centering
    \includegraphics{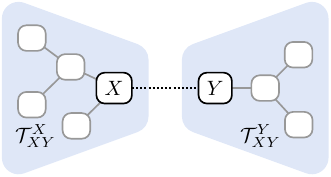}
    \caption{A tree decomposition $\td$ with bags $X$ and $Y$ and the tree decompositions $ \td^X_{XY} $ and $ \td^Y_{XY}$ highlighted in blue.
    }
    \label{fig:definition_T_XY}
\end{figure}

As an intermediate step to construct cycle-inducing tree decompositions, we first examine what happens if we split a tree decomposition into two parts by removing an edge.
For this let $\td$ be a tree decomposition, $XY \in E(\td)$, and $\td_1$ and $\td_2$ be the two components of $\td-XY$ such that $X \in V(\td_1)$. Then, we define $\td_{XY}^X = \td_1$ and $\td_{XY}^Y = \td_2$.
This is depicted in \cref{fig:definition_T_XY}.
Note that the following lemma is very similar to \cref{cor:inner_outer}.
In fact, \cref{cor:inner_outer} is a stronger statement for the special case that the tree decomposition is cycle-inducing.

\begin{lemma}\label{le:everything_on_one_side}
    Let $G$ be a plane triangulation with a tree decomposition $\td$ and let $ XY \in E(\td) $.
    Then, each face of $ G[X \cap Y] $ contains vertices of at most one of $ G[\td^X_{XY}] $ and $ G[\td^Y_{XY}] $ in its interior.
\end{lemma}

\begin{proof}
    First, assume that $S = X \cap Y$ is not a separator of $G$.
    Then, by the third condition of tree decompositions, one of the subtrees $\td^X_{XY}$ or $\td^Y_{XY}$ contains only vertices in $X \cap Y$.
    In particular, it contains no vertex in the interior of any face of $G[X\cap Y]$ and, thus, the statement holds.
    Therefore, we may assume that $S$ is a separator in $G$.
    Observe that any two vertices of $ G $ that are in the interior of the same face of $ G[S] $ are not separated by $ S $ since $ G $ is a triangulation.
    As $G[\td_{XY}^X]-S$ and $G[\td_{XY}^Y]-S$ are separated by $S$, each face contains vertices of at most one of them in its interior.
\end{proof}

Let $G$ be a graph and $\td$ be a tree decomposition of $G$ of width $k$.
We introduce the following notation which we use to proof \cref{thm:nice_edges}.
\begin{itemize}
    \item For a set $U \subseteq V(G)$ we define the \emph{deficiency} $d(U)$ of $U$ as the difference between the size of $ U $ and the smallest separator of $ G $ contained in $ U $.
    That is, minimal separators have deficiency $0$.
    If $U$ does not contain a separator of $G$, we set $d(U)= |U|$.
    \item For a natural number $i \in \mathbb{N}$ we define $d_i(\td)$ to be the number of different sets $U \subseteq V(G)$ with $d(U)=i$ such that there is an edge $XY \in E(\td)$ with adhesion set $U$, i.\,e., $X \cap Y = U$. Note that this is not the same as the number of edges $XY \in E(\td)$ such that $d(X \cap Y) = i$ as there might be a different edge $X'Y' \in E(\td)$ with the same adhesion set as $XY$. However, the number of such edges is an upper bound to $ d_i(\td) $.
    \item Finally, we define the \emph{deficiency sequence} $D(\td)$ of $\td$ as the tuple $(d_{k+1}(\td),\dots,d_1(\td))$ indicating how often each deficiency occurs. Note that there is no deficiency larger than $k+1$ as the bags of $\td$ have size at most $k+1$.
\end{itemize}

Note that if $ d(X \cap Y) = 0 $ for all edges $ XY $ in a tree decomposition of a triangulation, then each edge induces a minimal separator and, thus, a cycle by \cref{le:inclusion_minimal_separators_in_triangulations}.
That is, to find a cycle-inducing tree decomposition, we aim for a tree decomposition $\td$ with $ D(\td) = (0, \dots, 0) $.
To this end, we look in the following proof at a tree decomposition that is lexicographically minimal with respect to its deficiency sequence. 
A similar proof technique has been used before, e.g., by Bellenbaum and Diestel~\cite{induction_on_tree_decomposition}.

\cycleinducing*

\begin{proof}
    Let $\td$ be a tree decomposition of $G$ of width $k$ that is lexicographically minimal with respect to $D(\td)$, and let $XY \in E(\td)$ be an edge that maximizes $d(X \cap Y)$.
    If we have that $d(X \cap Y)=0$, we are done.
    So, we assume for the sake of contradiction that $d(X \cap Y) > 0$, and construct a tree decomposition $ \td'' $ of $ G $ that contradicts the minimality of $\td$.

    Before doing so, let us first state the following invariant, which we maintain throughout the proof for every tree decomposition $ \td^* $ we construct:
    Every edge $ A^* B^* \in E(\td^*) $ has a corresponding edge $ AB \in E(\td) $ whose endpoints are super-sets of $ A^* $, respectively $ B^* $. 
    We remark that multiple edges of $ \td^* $ may correspond to the same edge in $\td$.
    Note that this implies $ A^* \cap B^* \subseteq A \cap B $, so the difference to a minimal separator may only decrease compared to $\td$, i.\,e., we have $ d(A^* \cap B^*) \leq d(A \cap B) $.
    That is, having the invariant, it suffices to show that $ X \cap Y $ does not contribute to $ D(\td^*) $ for some tree decomposition $ \td^* $ to obtain $ D(\td^*) < D(\td) $.
    This then gives the desired contradiction to the minimality of $\td$.

    \begin{figure}
        \centering
        \includegraphics{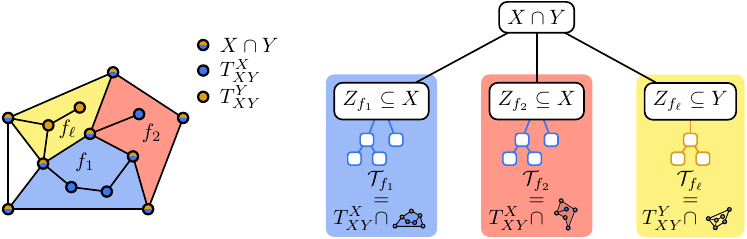}
        \caption{Left: The graph $G$ with chosen bags $X$ and $Y$ and faces $F=\{f_1,\dots,f_\ell\}$ of $G[X \cap Y]$. The white faces (including the outer face) are empty and, thus, can be omitted in the construction of $\td'$. Right: The augmented tree decomposition $\td'$.}
        \label{fig:cycle_inducing_tree_decomposition}
    \end{figure}

    As an intermediate step, we construct a tree decomposition $ \td' $ of $ G $.
    For this, we denote the set of faces of $G[X \cap Y]$ by $ F $ and do the following for each face $ f \in F $ as shown in \cref{fig:cycle_inducing_tree_decomposition}.
    First, we know by \cref{le:everything_on_one_side} that at most one of $\td_1 = \td_{XY}^X $ and $ \td_2 = \td_{XY}^Y $ has vertices inside $ f $.
    With this, let $\td_i$ with $i \in \{1,2\}$ be the tree decomposition that contains a vertex inside $f$ (if neither $\td_1$ nor $\td_2$ contains such a vertex, $f$ is a face of $G$ and we can skip $f$).
    Furthermore, let $\td_f$ be the tree decomposition resulting from intersecting every bag of $\td_i$ with $\interior_G[f]$ and let $Z_f $ be the bag of $\td_i$ that results from $X$, respectively $Y$.
    That is, we have $ Z_f = X \cap \interior_G[f] $, respectively $ Z_f = Y \cap \interior_G[f] $.
    Now, $\td'$ is obtained by taking the union of all $\td_f$ and connecting each $Z_f$ to a new bag $X \cap Y$.
    Note that this indeed yields a tree decomposition of $ G $: First, every vertex in $V(G) - (X \cap Y)$ is in exactly one face of $ G[X \cap Y] $ and, thus, every such vertex occurs in exactly one $\td_f$.
    In addition, every edge $xy \in E(G)$ is either in $ G[X \cap Y] $ or in some $f \in F$.
    In the first case $x$ and $y$ are contained in the new bag $X \cap Y$ and in the second case there is a bag in $T_f$ that contains $x$ and $y$.
    Furthermore, observe that the invariant indeed holds for $ \td^* = \td' $, since for the edges incident to the bag $ X \cap Y $ we may choose $ A = X $ and $ B = Y $, and the other edges have corresponding edges in $\td$ by construction.
    
    Next, we aim to get rid of edges $ CD \in E(\td') $ with $ C \subseteq D $.
    As long as there is such an edge, contract the smaller bag into the larger, and let $ \td'' $ denote the resulting tree decomposition.
    Note that $ \td'' $ indeed is a tree decomposition of the same width as $\td$ and our invariant is maintained.
    That is, we are left to show that $ X \cap Y $ does not contribute to $ D(\td'') $, i.\,e., that there is no edge whose adhesion set is $ X \cap Y $.
    For this, let $ Z \supseteq X \cap Y $ be the bag of $ \td'' $ corresponding to the bag $ X \cap Y $ of $ \td' $.
    We aim to show that $ Z $ is the only bag containing all vertices of $ X \cap Y $ as this implies that no edge has $ X \cap Y $ as its adhesion set.
    If $ Z = X \cap Y $, i.\,e., if no bag $ Z_{f} $ of $ \td' $ with $ f \in F $ is contracted into the bag $ X \cap Y $, this is evident from the construction of $ \td'' $.

    \begin{figure}
        \centering
        \includegraphics[page=2]{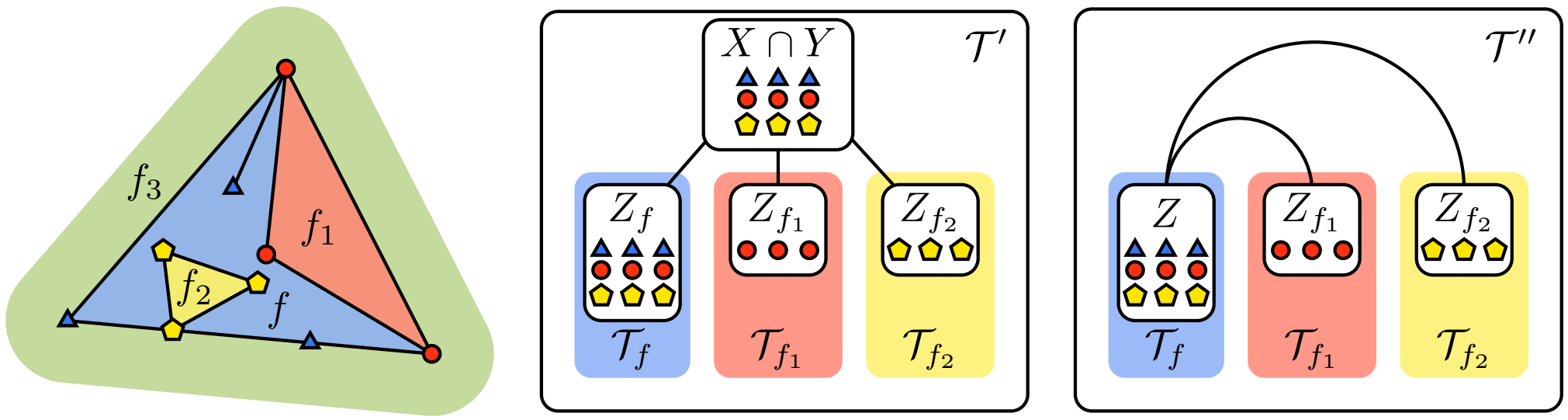}
        \caption{%
            Left: $ G[X \cap Y] $ with face $ f $ (blue) in the case that  $ Z $ results from contracting $ X \cap Y $ into some $ Z_f $.
            For ease of presentation we omit $f_3$ in $\td'$ and $\td''$.
            Middle: $ \td' $ with $ X \cap Y \subseteq Z_f $.
            Only the vertices in $X \cap Y$ are depicted in the bags.
            Right: $ \td'' $ after contracting $ X \cap Y $ into $ Z_f $.
        }
        \label{fig:cycle_inducing_tree_decomposition_outerplanar}
    \end{figure}

    So assume that $ Z $ results from contracting $ X \cap Y $ into some $ Z_f $ of $ \td' $ (and possibly further contractions), where $ f \in F $ is some face of $ G[X \cap Y] $.
    In particular, we have $ X \cap Y \subseteq Z_f $.
    Since $ Z_f $ is defined as $ X \cap \interior_G[f] $, respectively $ Y \cap \interior_G[f] $, this means that all vertices of $ X \cap Y $ are contained in $ \interior_G[f] $.
    As $ f $ is a face of $ G[X \cap Y] $, we conclude that $ G[X \cap Y] $ is outerplanar and $ f $ is incident to all its vertices.
    We refer to \cref{fig:cycle_inducing_tree_decomposition_outerplanar} for the current situation.
    Recalling our assumption that $ d(X \cap Y) > 0 $ and, thus, that $G[X \cap Y]$ is not an induced cycle by \cref{le:inclusion_minimal_separators_in_triangulations}, we obtain that $ f $ is the only such face in $ F $.
    In particular, no other face $ f' \neq f $ of $ G[X \cap Y] $, and therefore no other $ Z_{f'} $, contains all vertices of $ X \cap Y $.
         
    Next, we apply \cref{le:everything_on_one_side} to show that $ Z $ does not have a neighbor $ Z' $ with $ X \cap Y \subseteq Z' $, which is what we aimed for.
    Since $ f $ is the only face of $ G[X \cap Y] $ containing all vertices of $ X \cap Y $, we only need to search for $ Z' $ in the bags resulting from $ \td_f $.
    To do so, assume there is such a neighbor $Z'$, and consider the edge $ Z Z' $ in $ \td'' $.
    We start by proving that the adhesion sets of $ XY $ and of $ Z Z' $ are the same, i.\,e., that $ X \cap Y = Z \cap Z' $.
    First, we have $ X \cap Y \subseteq Z \cap Z' $ by definition of $ Z $ and $ Z' $.
    Second, by our invariant, there is some edge $ AB $ in $\td$ with $ Z \cap Z' \subseteq A \cap B $.
    Now, if $ X \cap Y \subsetneq Z \cap Z' $, then the adhesion set of $ AB $ is a strict super-set of $ X \cap Y $, which contradicts the maximality of $ d(X \cap Y) $.
    Hence, we have $ X \cap Y = Z \cap Z' $.

    Having this, we apply \cref{le:everything_on_one_side} to $ \td'' $ and the edge $ Z Z' $.
    We obtain that either $ Z $ or $ Z' $ do not contain any vertex in the strict interior of any face of $ G[Z \cap Z'] = G[X \cap Y] $, and in particular of $ f $. 
    However, recall that all vertices in $Z$ and $Z'$ are in $\td_f$, i.\,e., they are incident to $f$ or in the strict interior of $f$.
    Thus, either $Z$ or $Z'$ are a subset of $Z \cap Z'$. That is, we have $ Z \subseteq Z' $ or $ Z' \subseteq Z$, contradicting the construction of $ \td'' $ since the edge $ Z Z' $ can be contracted.

    To sum up, we get that $ Z $ does not have a neighbor $ Z' $ containing $ X \cap Y $ and, thus, $ \td'' $ contains no edge with adhesion set $X \cap Y $.
    Together with our invariant, this shows $ D(\td'') < D(\td) $, contradicting the minimality of $\td$.
    We conclude that $ D(\td) = 0 $ and, hence, $\td$ is a cycle-inducing tree decomposition.
\end{proof}

\section{Non-Strong Planar Quasi-\texorpdfstring{$\bm{k}$}{k}-Trees} \label{sec:niceness}

Recall that a planar quasi-$k$-tree is a \nice planar quasi-$k$-tree if the graph is biconnected at every step of the construction (cf.\ \cref{def:strong_quasi_k-trees}).
Complementing \cref{thm:nice_4-trees}, we aim to show that for every $k\geq 5$ there is a planar graph with treewidth $5$ that is not a subgraph of a \nice planar quasi-$k$-tree.

We emphasize two differences of \cref{def:strong_quasi_k-trees} to \cref{def:quasi_k-trees}:
First, the base case is always a triangle instead of a cycle of length up to $ k $ as the triangle is the only cycle that contains a biconnected subgraph on three vertices.
Second and most importantly, we maintain biconnectivity in every step of the construction.
We remark that starting with a triangle is not a restriction as it is easy to build a larger cycle as a subgraph by repeatedly inserting vertices of degree 2.
That is, the crucial difference is the biconnectivity, and indeed the graphs we construct for \cref{thm:not_nice_k_tree} to separate planar quasi-$k$-trees from strong ones are triangulations, i.\,e., start with a triangle.

We first describe graphs $ H_k $ and $ G_k $ based on which we construct a triangulated graph $G'_k$ that is not a strong planar quasi-$k$-tree.
Then, at the end of the section, we show that there is not even a supergraph that is a strong planar quasi-$k$-tree.
We do so by showing in \cref{le:supergraphs} that a triangulation is a subgraph of a \nice planar quasi-$k$-tree if and only if it is a \nice planar quasi-$k$-tree itself. 
Since the graph we construct is a triangulation, \cref{thm:not_nice_k_tree} then follows immediately.

\begin{figure}
    \centering
    \includegraphics{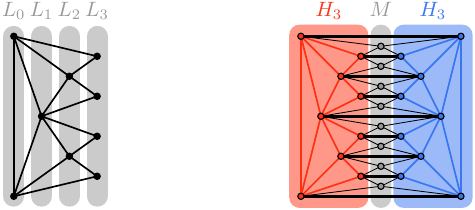}
    \caption{Left: The graph $H_3$ with layers $L_0,\dots,L_3$ highlighted in gray. The graph induced by the layers $L_0,\dots,L_i$ is the graph $H_i$.
    Right: The graph $G_3$ with layer $M$ highlighted in gray. 
    The two copies of $H_3$ are depicted in red and blue, respectively, and the edges connecting the copied vertices are thick black.}
    \label{fig:outerplanar}
\end{figure}

We inductively define the graph $H_k$ as follows, see \cref{fig:outerplanar}~(left).
The graph $H_0$ is $K_2$, i.\,e., a single edge.
We say all vertices of $H_0$ are on layer $L_0$.
To obtain $ H_k $ for $ k > 0 $, take $ H_{k - 1} $ and for every edge with an endpoint in the $(k-1)$-st layer, add a vertex that is connected to both endpoints of the edge.
We refer to these newly added vertices as the $k$-th layer $L_k$ of $H_k$.
All other vertices inherit the layers from $H_{k-1}$.
Note that every $H_k$ is outerplanar.
We refer to the edge of $H_k$ that is induced by $L_0$ as the \emph{base edge} of $H_k$.

To construct $ G_k $, we take two disjoint copies of $H_k$ and connect each pair of vertices that are copies of each other, see \cref{fig:outerplanar}~(right).
Now, we fix a planar embedding of this graph such that the two base edges of the $H_k$-copies are incident to the outer face.
Observe that each inner face is either a triangle induced by a single copy of $ H_k $ or a quadrangle bounded by some edge of $ H_k $, its copy, and two edges connecting the copies.
Finally, for every inner quadrangle, we add a vertex that is adjacent to all vertices incident to the face.
We refer to the set of these vertices as layer $M$ of $G_k$, and to the two base edges of the $H_k$-copies as the \emph{base edges} of $G_k$.

\begin{figure}
    \centering
    \includegraphics{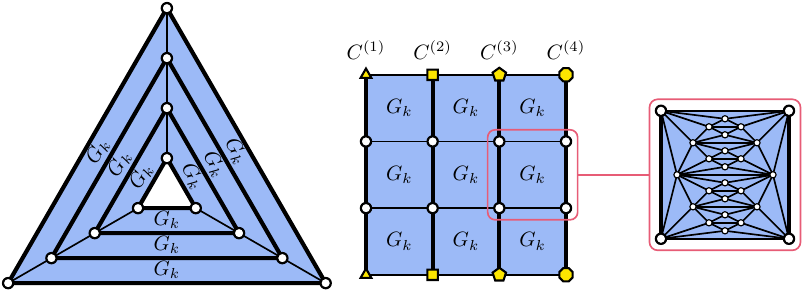}
    \caption{Left and Middle: Two different embeddings of the graph $G_k'$. Each yellow vertex is identified with the other yellow vertex of the same shape. Each quadrangular face of $C_3 \Box P_4$ (blue) is identified with a copy of $G_k$. The edges of the copies $C^{(1)},\dots,C^{(4)}$ of the cycle $C_3$ are thick black.
    Right: A closeup of a copy of $G_k$ for $k=3$. The base edges of $G_k$ that are identified with the edges of copies of $C_3$ are thick black.}
    \label{fig:not_nice_k_tree}
\end{figure}

Finally, we construct $G_k'$.
For this, consider a path $P_4 = (v_1,v_2,v_3,v_4)$ and the Cartesian product $C_3\Box P_4$ as shown in \cref{fig:not_nice_k_tree}.
This graph is $3$-connected and, thus, has a unique combinatorial embedding.
For a vertex $v_i \in V(P_4)$, we refer to the graph induced by the set $\{(v_i,u) \colon u \in V(C_3)\}$ as a \emph{copy} of $C_3$ and denote it by $C^{(i)}$.
The copies of $C_3$ are thick black in \cref{fig:not_nice_k_tree}.
For each quadrangular face of $C_3\Box P_4$, we introduce a copy of $ G_k $ and identify the four endpoints of the base edges with the four vertices of the face. 
We do this so that the base edges of $G_{k}$ are identified with edges of copies of $C_3$ (this means the base edges are vertical in the middle of \cref{fig:not_nice_k_tree}).
Note that $G_k'$ is a triangulation as $C_3\Box P_4$ only has triangular and quadrangular faces and $G_k$ is internally triangulated.

Before showing that $ G_k' $ is not a strong planar quasi-$k$-tree, we bound its treewidth.

\begin{restatable}{lemma}{twnotnice}\label{le:G_i_small_tw}
    For every $ k \geq 0 $, the graph $G_k'$ has treewidth at most $5$.
\end{restatable}

\begin{proof}
    To start, we show that the graph $G_k$ has treewidth at most $5$ and admits a certain tree decomposition.
    Then, the treewidth bound of $G_k'$ follows.

    \begin{claim}\label{claim:nice_tree_decomposition}
        For every $ k \geq 0 $, the graph $G_k$ has a tree decomposition of width $5$ that contains a bag with the four vertices incident to the base edges of $G_k$.
    \end{claim}
    
    \begin{claimproof}
        As $H_k$ is outerplanar, it has treewidth at most $2$ and, thus, a tree decomposition with bag size at most $3$.
        We take such a tree decomposition $\td$ of $H_k$ and augment it to obtain a tree decomposition $\td''$ of $G_k$.
        For each vertex $v$ in a bag of $\td$, we substitute $v$ with the two copies of $v$ in $G_k$.
        This gives us a valid tree decomposition $ \td' $ of $G_k$ except for the vertices in layer $M$ (depicted gray in \cref{fig:outerplanar} on the right), where each bag contains at most six vertices.
        To construct $ \td'' $, recall that each vertex $ w \in M$ is only adjacent to four vertices resulting from two adjacent vertices $ x, y $ of $ H_k $.
        As there is a bag in $\td$ that contains the original edge $ xy $, there is a bag $ B $ in $\td'$ that contains both copies of $ xy $.
        Now, we add a new bag adjacent to $ B $ that consists of $ w $ and the two copies each of $ x $ and $ y $, so five vertices in total.
        Doing so for each vertex of $ M $ yields the desired tree decomposition $ \td'' $ of $ G_k $ with bag size at most 6.
        As there is a bag in $\td$ that contains the base edge of $H_k$, there is a bag in $\td'$, and, thus, also in $ \td'' $, that contains both of its copies and, hence, the base edges of $G_k$.
    \end{claimproof}
    \begin{figure}
        \centering
        \includegraphics{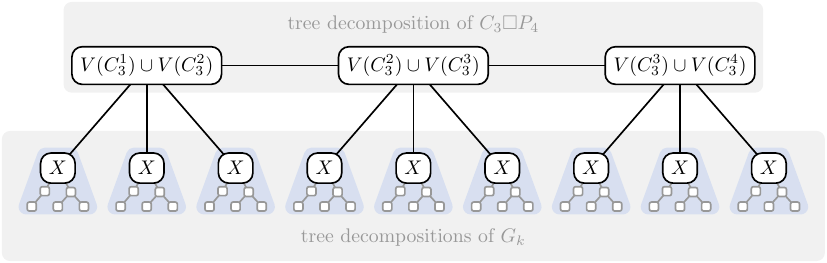}
        \caption{The tree decomposition of the graph $G_k'$. The path on the top is a tree decomposition of $C_3 \Box P_4$. The bags denoted with $X$ are the bags containing the base edges of the $G_k$ copies. The existence of these bags is guaranteed by \cref{claim:nice_tree_decomposition}.}
        \label{fig:tree_decomposition_G}
    \end{figure}
    To bound the treewidth of $G_k'$ recall that $C^{(1)},\dots,C^{(4)}$ denote the copies of $C_3$ in $C_3 \Box P_4$.
    We start with a path decomposition of $C_3 \Box P_4$ with bags $X_i = V(C^{(i)}) \cup V(C^{(i+1)})$ such that each $X_i$ is a neighbor of $X_{i+1}$ for $1 \leq i \leq 3$.
    As each $C^{(i)}$ contains three vertices, the path decomposition has width $5$.
    Furthermore, for every quadrangular face of $C_3 \Box P_4$ there is a bag containing its incident vertices.
    By \cref{claim:nice_tree_decomposition}, each copy of $G_k$ has a tree decomposition of width $5$ that contains the four vertices that are endpoints of the base edges of $G_k$ in a single bag.
    Making this bag adjacent to the appropriate bag in the tree decomposition of $C_3 \Box P_4$ yields a tree decomposition of $G$ of width $5$. For a depiction see \cref{fig:tree_decomposition_G}.
\end{proof}

To show that $G_k'$ is not a \nice planar quasi-$k$-tree, we look at vertex-orderings of subgraphs $H$ of $H_k$.
Think about these orderings as potential \nice construction sequences of $H$.
We show that if $H$ contains vertices from both the first layer $L_0$ and the last layer $L_k$, in order for such a \nice construction sequence to exist, $H$ also contains vertices from all other layers $L_1,\dots,L_{k-1}$.
This lemma is crucial for showing that $G_k'$ is not a strong planar quasi-$k$-tree, as we show that in every potential strong $k$-construction sequence of $G_k'$, such a graph $ H $ occurs in some step, so that the vertices we obtain from the lemma lie on a common face of size $ k + 1 $.
Thus, the sequence is not a $ k $-construction sequence.

\begin{restatable}{lemma}{largeface}\label{le:large_face}
    Let $ k \geq 0 $ and let $\sigma = (v_1,\dots,v_n)$ be a vertex-ordering of a subgraph $H$ of $H_k$, where $ n = |V(H)| $, such that 
    \begin{enumerate}
        \item $\sigma$ starts with the two vertices in $L_0$ of $H_k$, 
        \item the graph $H[\{v_1,\dots,v_i\}]$ is biconnected for every $3\leq i \leq n$, and
        \item $H$ contains a vertex that is in layer $L_k$ of $H_k$.
    \end{enumerate}
    Then, there are vertices $u_0,\dots, u_k$ with $u_j \in L_j $ for $ j = 0, \dots, k $ such that $u_0 <_\sigma \dots <_\sigma u_k$.
\end{restatable}

\begin{proof}
    For a vertex $u \in L_j$ with $ 1 \leq j \leq k$, we call the neighbors of $u$ in $L_{j-1}$ the \emph{parents} of $u$.
    We show that for each vertex $ u \in V(H) \setminus L_0 $, its parents are also in $ H $ and occur before $ u $ in $ \sigma $.
    For this, let $ H' \subseteq H $ be the subgraph of $ H $ induced by all vertices until $ u $ according to $ \sigma $.
    By the first two properties, $ H' $ is biconnected.
    Note that the two parents of $u$ separate $u$ from the vertices in $L_0$ in $H'$.
    Thus, in order for $H'$ to be biconnected, both parents of $u$ are in $H'$.
    By definition of $ H' $, we indeed obtain that the parents of $ u $ occur before $u$ in $\sigma$.
    
    To finish the proof, recall that by definition of $ H_k $, each vertex has a parent in the previous layer.
    By the third property, $ H $ contains a vertex $ u_k \in L_k $.
    We choose $ u_{k-1}, \dots, u_0 $ by taking a parent from the previous layer in each step, where the desired order with respect to $ \sigma  $ is ensured by our claim above.
\end{proof}

Now, we come to the main part of the proof of \cref{thm:not_nice_k_tree}.

\begin{figure}
    \centering
    \includegraphics{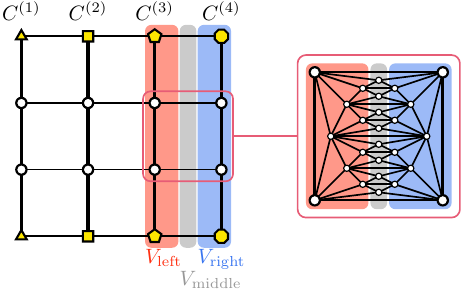}
    \caption{Left: The graph $G_k'$. 
    The sets $V_\text{left},V_\text{middle}$, and $V_\text{right}$ are highlighted in red, gray and blue, respectively.
    Right: A closeup of a $(3,4)$-copy of $G_k$ for $k=3$. The base edges of $G_k$ that are identified with the edges of $C_3$ copies are thick black.}
    \label{fig:finding_a_large_face}
\end{figure}

\begin{lemma}\label{le:not_nice}
    For every $k \geq 5$, the graph $G_k'$ is not a \nice planar quasi-$k$-tree.
\end{lemma}

\begin{proof}
    Recall that $ G_k' $ is constructed from $ C_3 \Box P_4 $ by inserting copies of $ G_k $ into each quadrangular face.
    For $ i = 1, 2, 3 $, we call a copy of $ G_k $ whose base edges are identified with edges from $ C^{(i)} $ and $ C^{(i+1)} $ a \emph{$ (i, i+1) $-copy}, where $ C^{(1)}, \dots, C^{(4)} $ denote the four copies of $ C_3 $ in the product.
    Assume for the sake of contradiction that $G_k'$ is a \nice planar quasi-$k$-tree.
    Then, there is a \nice $k$-construction sequence $\sigma = (v_1,\dots, v_n)$ of $G_k'$.
    Without loss of generality let $v_1$ not be in a $(3,4)$-copy of $G_k$,
    i.\,e., $v_1$ is to the left of $C^{(3)}$ in \cref{fig:not_nice_k_tree}.
    Using the biconnectivity of every step of the $k$-construction sequence, we aim to find a face of size $ k + 1 $ in some step, contradicting our assumption as this is not allowed in a $k$-construction sequence. 
    We do this by finding a suitable subgraph $H$ of some copy of $H_k$ for which we can apply \cref{le:large_face}.
    Then, the $ k $ vertices provided by the lemma together with a vertex from layer $ M $ of the respective copy of $ G_k $ form the desired face of size $ k + 1 $.

    As shown in \cref{fig:finding_a_large_face}, we partition the vertices of the $ (3,4) $-copies into three sets based on the two copies of $ H_k $ and the additional layer $ M $ each $ G_k $ consists of.
    The sets $ V_\text{left} $, respectively $ V_\text{right} $, contain the vertices of the $ H_k $-copies whose base edge is in $ C^{(3)} $, respectively in $ C^{(4)} $, and the set $ V_\text{middle} $ contains the remaining vertices, i.\,e., those that are in layer $ M $ of some $ G_k $.
    As an intermediate step, we show that the first vertex $ m_1 $ of $ V_\text{middle} $ occurs before the first vertex of $ V_\text{right} $ in $ \sigma $, and then find $ H $ among vertices of $ V_\text{left} $ that are inserted before $ m_1 $.
    For this, first note that to insert the first vertex of $ V_\text{right} $ with insertion degree 2, two of its neighbors need to be inserted beforehand.
    However, we show that there is only one candidate that is neither in $ V_\text{right} $ nor in $ V_\text{middle} $.
    Indeed, the vertices of $ V_\text{right} $ have only neighbors in $ (3,4) $-copies of $ G_k $, and by construction of $ G_k $, all but at most one are in the same copy of $ H_k $ or in $ M $, i.\,e., also in $ V_\text{right} \cup V_\text{middle} $.
    Hence, in order to insert the first vertex of $ V_\text{right} $ in the construction sequence, we need at least one of its neighbors in $ V_\text{middle} $ to occur first.
    That is, $ m_1 $ is the first vertex in $ V_\text{middle} \cup V_\text{right} $.
    Moreover, by construction, $ m_1 $ has exactly two neighbors that are not in $ V_\text{middle} \cup V_\text{right} $, both in the same copy of $ H_k $, so these two neighbors occur before $ m_1 $ in $ \sigma $ for $ m_1 $ to have insertion degree 2.
    Let $ H_k' $ denote the copy of $ H_k $ containing these two neighbors.

    \begin{figure}
        \centering
        \includegraphics{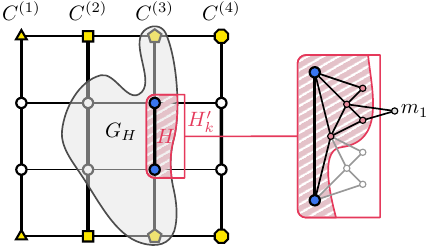}
        \caption{Left: The graph $G_H$ (gray), the graph $H$ (hatched red), and the graph $H_k'$ (red).
        The vertices incident to the base edge of $H_k'$ are highlighted in blue.
        Right: A closeup of $H_k'$.
        The vertex $m_1$ is depicted in gray.}
        \label{fig:H_biconnected}
    \end{figure}

    Now, heading for an application of \cref{le:large_face}, let $ H $ be the subgraph of $ H_k' $ induced by the vertices before $ m_1 $ in $ \sigma $, and let $ G_H $ be the subgraph of $ G_k' $ before $ m_1 $ is inserted.
    We refer to \cref{fig:H_biconnected} for an example.
    In particular, we have $ H \subseteq G_H $, neither of them contains any vertices of $ V_\text{middle} $ or $ V_\text{right} $, and $ G_H $ is biconnected as it is built according to a strong $k$-construction sequence.
    It follows that the the base edge of $ H_k' $ is contained in $ G_H $ and induces a separator of size 2.
    
    Using this, we argue that $ H $ indeed satisfies the requirements of \cref{le:large_face}.
    For the first condition, we need to show that the endpoints of the base edge of $ H_k' $ are inserted before any other vertex of $ H $ in $ \sigma $.
    For this, recall that $v_1$, i.\,e., the very first vertex of $ \sigma $, is not in a $(3,4)$-copy of $G_k$ and, thus, not in $ H $.
    That is, removing the separator induced by the base edge of $ H_k' $ leaves two components, one containing $ v_1 $ and the other consisting of $ H $ except for the base edge.
    Hence, since $ G_H $ is biconnected at every step of the construction, the base edge is inserted before any other vertex of $ H $.
    The second condition asks for biconnectivity at every step of the construction sequence.
    We already have this property for $ G_H $ and now need to transfer it to $ H $.
    Observe that the base edge of $ H_k' $ is a separator of size 2 not only in $ G_H $ but in each construction step after inserting the base edge until $ G_H $ is finished.
    Hence, in each such step, each component plus the separator is biconnected.
    In particular, the component induced by the vertices of $ H $ that are already inserted is biconnected, as required.
    For the third condition, recall that $ H_k' $ is defined so that $ m_1 $ has two neighbors in $ H_k' $, both before $ m_1 $ in $ \sigma $ and thus in $ H $.
    Further recall that $ m_1 $ is in layer $ M $ of the respective copy of $ G_k $, so by construction of $G_k$, one of the neighbors is in layer $L_k$ of $H_k'$.

    To sum up, we have all three conditions for \cref{le:large_face} and obtain vertices $u_1,\dots,u_k$ in $H$.
    As $ H_k' $, and thus also $ H $, is outerplane in the embedding of $G_k'$, we have that $u_1,\dots,u_k$ share a face in $ G_H $ into which $ m_1 $ is inserted in the next step.
    However, $ m_1 $ has exactly two neighbors in $ H $ that are consecutive on the boundary of the face that $m_1$ is inserted into.
    Hence, the size of the face is increased to at least $ k + 1 $ by $ m_1 $.
    This contradicts $\sigma$ being a $k$-construction sequence.
\end{proof}

Now, as the last step to prove \cref{thm:not_nice_k_tree}, we show that a planar triangulation is a subgraph of a \nice planar quasi-$k$-tree if and only if it is a \nice planar quasi-$k$-tree itself.
This gives us that $G_k'$ is not a subgraph of a \nice planar quasi-$k$-tree and, thus, \cref{le:not_nice,le:supergraphs} together prove \cref{thm:not_nice_k_tree}.

\begin{restatable}{lemma}{supergraphs}\label{le:supergraphs}
    Let $G$ be a triangulation.
    Then, there is a supergraph of $G$ that is a \nice planar quasi-$k$-tree if and only if $G$ is a \nice planar quasi-$k$-tree.
\end{restatable}

\begin{proof}
    Let $H$ be a supergraph of $G$ that is a \nice planar quasi-$k$-tree and $\sigma =(v_1,\dots,v_n)$ be a \nice $k$-construction sequence of $H$.
    We fix a canonical embedding of $H$.
    Let $G^\sigma_i$ be the maximal induced subgraph of $H^\sigma_i$ that contains only vertices of $G$.
    Our goal is to show that each $G^\sigma_i$ with at least three vertices is biconnected and has faces of size at most $k$.
    It then follows, that $G$ is a \nice planar quasi-$k$-tree.

    Consider any $G^\sigma_i$ with at least three vertices and assume that it is not biconnected.
    Then, there is a cut-vertex $c$ such that $G^\sigma_i$ has at least two components.
    As $H^\sigma_i$ is biconnected, vertex $c$ is not a cut-vertex in $H^\sigma_i$.
    Let  $A$ and $B$ be two components of $G^\sigma_i$ with minimal distance in $H^\sigma_i$.
    By the choice of $A$ and $B$, there is a path connecting $A$ and $B$ in $H^\sigma_i - c$ that is disjoint from $G_i^\sigma$ except for its end-vertices.
    Due to planarity, this path is embedded inside a single face of $G^\sigma_i$ and $G$.
    Thus, there is a face in $G$ that is incident to $A$ and $B$.
    As $G$ is a triangulation $A$ and $B$ are adjacent.
    This contradicts $c$ being a cut-vertex in $G^\sigma_i$.
    Thus, $G^\sigma_i$ is biconnected.

    Now let $f$ be a face of $G^\sigma_i$.
    If $f$ is a triangle we are done.
    Consider the subgraph $H_f$ of $H^\sigma_i$ that is embedded inside $f$.
    Due to planarity each component of $H_f$ is embedded inside a single face of $G$.
    As each face of $G$ is a triangle but $f$ is not a triangle, each component of $H_f$ is adjacent to at most two vertices on the boundary of $f$.
    Moreover, these vertices are adjacent on the boundary of $f$.
    Thus, there is a face $f'$ in $H^\sigma_i$ that is incident to all vertices (and possibly more) that $f$ is incident to.
    Since $f'$ has size at most $k$, so does $f$.

    Together, this gives us that $G$ is a \nice planar quasi-$k$-tree.
    For the other direction we can take $G$ itself as the supergraph and are done.
\end{proof}

\begin{remark}
    While \cref{thm:not_nice_k_tree} states that there is a graph of treewidth $5$ that is not a subgraph of a \nice planar quasi-$k$-tree, one can easily obtain such a graph of treewidth $t$ for all $t \geq 5$ by taking the disjoint union of the graph constructed in the proof of \cref{le:not_nice} and any planar graph of treewidth $t$ and triangulating it without increasing the treewidth.
\end{remark}

\section{Open Questions}\label{se:open_questions}

By \cref{thm:quasi_k_trees}, for every $k \geq 3$, a planar graph~$G$ has treewidth at most~$k$ if and only if $G$ is a subgraph of a planar quasi-$k$-tree.
By \cref{thm:nice_4-trees,thm:not_nice_k_tree}, this is equivalent to $G$ being a subgraph of a \nice planar quasi-$k$-tree if and only if $k \in \{3,4\}$.
It may be worthwhile to investigate the subgraphs of \nice planar quasi-$k$-trees for $k > 4$.

\begin{question}
    For $k > 4$, which graphs are subgraphs of \nice planar quasi-$k$-trees?
\end{question}

For several important subclasses of planar graphs, including planar $3$-trees, the program \texttt{plantri}~\cite{plantri3,plantri1,plantri2} by Brinkmann and McKay efficiently generates all graphs of that class.
In order to avoid costly isomorphism checking, the program uses inductive definitions of the classes.
With \cref{def:quasi_k-trees} we give for every $k \geq 4$ the first inductive definition of the class of all triangulations of treewidth at most~$k$.

\begin{question}
    Can planar quasi-$k$-trees be used to efficiently generate all planar triangulations of treewidth at most~$k$ on $n$ vertices?
\end{question}

It is also interesting how fast a $k$-construction sequence can be computed.
In the next paragraph we briefly discuss that, having a cycle-inducing tree decomposition, the proof of \cref{le:quasi_k-trees} can be implemented in $ O(k \cdot n) $.
That is, the crucial step is to compute a cycle-inducing tree decomposition.
Note that our proof for their existence is constructive (\cref{le:inclusion_minimal_separators_in_triangulations}), but not optimized for efficiency and we do not investigate how many steps the construction takes.
We expect that computing a $ k $-construction sequence is FPT in~$ k $.

\begin{question}
    For a planar triangulation of treewidth $ k $, how efficiently can a $k$-construction sequence be computed?
\end{question}

To compute a $ k $-construction sequence for a given triangulation $G$ with $n$ vertices, we assume to be given a cycle-inducing tree decomposition of width $k$ with a linear number of bags.
Note that such a tree decomposition exists by \cref{le:treewidth}, where the constructed tree decomposition indeed has a 1-to-1 correspondence between inserted vertices and bags.
To obtain the $k$-construction sequence, for each bag $B$ we get the cycle $(v_1,\dots,v_\ell)$ bounding $G[B]$ by computing the intersection $A \cap B$ with the parent $A$ of $B$.
Then, we use a BFS, starting with the vertices $v_1,\dots,v_\ell$, to obtain some ordering $\sigma_B=(v_1,\dots,v_{|B|})$ of $B$ such that each $G[\{v_1,\dots,v_i\}]$ is connected for $ 1 < i \leq |B| $.
We define $\sigma_B'= (v_{\ell+1},\dots,v_{|B|})$.
Finally, using a bottom-up approach, for each bag $B$ with children $C_1,\dots,C_d$, we append $\sigma_{C_1}'\!\cdot \ldots \cdot \sigma_{C_d}'$ to $\sigma_B'$.
Appending the resulting sequence of the root to any ordering of the outer triangle yields a $k$-construction sequence.
For each bag this takes time $O(k)$.
Since we have a linear number of bags, it follows that the algorithm takes time $O(k\cdot n)$.

Finally, we conclude the paper with a conjecture regarding the insertion degree of construction sequences.
It is well known that planar graphs have degeneracy at most $5$, i.\,e., for each planar graph there is a vertex-ordering such that every vertex~$v$ has at most five neighbors preceding~$v$ in this ordering.
Our definition of planar quasi-$k$-trees yields a $k$-construction sequence where every vertex has at most~$k$ preceding neighbors, i.e., each vertex has insertion degree at most~$k$.
We believe that we can show that every planar quasi-$k$-tree admits a $k$-construction sequence such that each vertex has insertion degree at most~$O(\sqrt{k})$.
However, our approach is quite tedious and seemingly can not be used to show any bound better than $O(\sqrt{k})$.
Since we do not think that this is asymptotically optimal, we do not think that our proof would give any relevant insight and omit it here.
Instead, we conjecture the following.

\begin{conjecture}\label{con:constant_insertion_degree}
    There is a universal constant $c\in \mathbb{N}$ such that for every planar quasi-$k$-tree there is a $k$-construction sequence such that every vertex has insertion degree at most~$c$.
\end{conjecture}
We believe that this constant is $5$.

\bibliography{references}

\appendix

\section[Comparing Planar Quasi-\texorpdfstring{$k$}{k}-Trees and Sphere-Cut-Decompositions]{Comparing Planar Quasi-\texorpdfstring{$\bm{k}$}{k}-Trees  and Sphere-Cut-Decompositions} \label{sec:comparing_sphere_cut}

In this section, we discuss the similarities and differences of $k$-construction sequences to sphere-cut-decompositions.
(For a definition of sphere-cut-decompositions we refer to~\cite{sphere_cut}.)
A result by Seymour and Thomas~\cite{ratcatcher} and Dorn, Penninkx, Bodlaender, and Fomin~\cite{sphere_cut} gives that every plane triangulation of branchwidth at most~$b$ has a sphere-cut-decomposition of width at most~$b$.
As we describe next, such a sphere-cut-decomposition in particular gives a construction sequence for the plane triangulation.

Let $G$ be a plane triangulation.
A \emph{noose} is a subset of the plane that is homeomorphic to a circle that intersects $G$ only in vertices of $G$.
The length of a noose is the number of vertices of $G$ it meets.
We start with a noose meeting all (three) outer vertices of $G$.
Then, we iteratively replace one noose by two nooses, while maintaining that every noose has length at most~$b$.
Namely, having a noose of length at most~$b$, we add a chord in its interior, connecting two points of the noose, such that each of the two resulting inner faces is bounded by a noose of length at most~$b$.
We do this, replacing one noose by two nooses in each step, until every vertex of $G$ is met by some noose.

To compare this with planar quasi-$k$-trees, think of this procedure with nooses as a construction sequence of $G$ in the following sense:
We start with the cycle consisting of the (outer) vertices met by the starting noose.
Then, every time we insert a new chord, we insert the vertices the chord meets.
Edges are added as soon as both incident vertices are present.
This is very similar to a $k$-construction sequence.
In fact, a sequence of nooses of length at most~$b$ immediately yields a $k$-construction sequence of $G$ with $k = \floor{\frac{3}{2}b}-1$.
This follows, as the length of a noose plus the length of a chord that is inserted into it has at most size $\floor{\frac{3}{2}b}$ (this is attained for a noose of length~$b$ and a chord of length~$\floor{\frac{b}{2}}$).
Conversely, every $k$-construction sequence of $G$ yields a sequence of nooses of length at most~$b = k+1$.

To conclude, the construction sequences obtained by sphere-cut-decompositions are quite similar to the $k$-construction sequences we describe.
However, they differ in some details.
Most importantly, the size of a noose plus a chord in the construction with nooses and the size of a face for the construction of treewidth differ by a factor of up to $\frac{3}{2}$.
We believe that this difference can be significant especially for graphs of small treewidth.
For example, it follows from \cref{thm:nice_4-trees} that planar graphs of treewidth~$4$ have branchwidth at most~$4$, as a \nice $4$-construction sequence can easily be transformed into a sequence of nooses of length at most~$4$ as described above.
However, in a general sequence of nooses of length at most~$4$, there might be a chord of length~$2$ inserted into a noose of length~$4$, increasing the number of \enquote{active} vertices at each step from at most~$4$ to at most~$5$ and, thus, increasing the number of cases drastically for proofs such as the proof of \cref{thm:arboricity}.
On the other hand, if a triangulation $G$ has branchwidth~$b$ and treewidth~$k$ and we have $k = \floor{\frac{3}{2}b}-1$, then a sphere-cut-decomposition of $G$ is basically equivalent to a $k$-construction sequence of $G$.
\end{document}